\documentclass[leqno]{siamart171218}

\usepackage[framemethod=TikZ]{mdframed}
\usepackage[final=true]{microtype}
\usepackage[strict]{changepage}
\usepackage[english]{babel}
\usepackage{hyperref}
\usepackage[all]{hypcap}
\usepackage{transparent}
\usepackage{mathtools}
\usepackage{adjustbox}
\usepackage{rotating}
\usepackage{multirow}
\usepackage{csquotes}
\usepackage{arydshln}
\usepackage{dashrule}
\usepackage{verbatim}
\usepackage{placeins}
\usepackage{stmaryrd}
\usepackage{graphicx}
\usepackage{wrapfig}
\usepackage{amssymb}
\usepackage{amsmath}
\usepackage{amsbsy}
\usepackage{xcolor}
\usepackage{ifthen}
\usepackage{framed}
\usepackage{eufrak}
\usepackage{color}
\usepackage{tikz}
\usepackage{bbm}
\usepackage{bm}
\usepackage{xspace}

\allowdisplaybreaks

\newcommand{\I}{{\textit{II}}}
\newcommand{\II}{{\textit{I}}}

\newcommand{\IIa}{{\textit{Ia}}}
\newcommand{\IIb}{{\textit{Ib}}}
\newcommand{\Edit}[1]{{#1}}

\newcommand{\EEdit}[1]{{#1}}

\renewcommand{\div}{\mathrm{div}} 
\newcommand{\Pk}{{\mathbb{P}_{k}}}
\newcommand{\Pkk}{{\mathbb{P}_{[k]}}}
\newcommand{\Qk}{\mathbb{Q}_{k}}
\newcommand{\Qkk}{\mathbb{Q}_{[k]}}
\newcommand{\R}{\mathbb{R}}
\newcommand{\dd}{\mathrm{d}}
\newcommand{\Rtt}{RT_{k}(K)}
\newcommand{\Rk}{\mathcal{R}_{k}(\partial K)}
\newcommand{\hk}{\mathcal{H}_{k}(\partial K)}
\newcommand{\Dofs}{degrees of freedom }
\newcommand{\Dofss}{degrees of freedom}
\newcommand{\Int}[4]{\int\limits_{#1}^{#2}#3\,\mathrm{d}#4}
\newcommand{\IntL}[4]{\int_{#1}^{#2}#3\,\mathrm{d}#4}

\newcommand{\x}{\mathrm{x}}
\newcommand{\y}{\mathrm{y}}

\newcommand{\Hk}{\mathbb{H}_{k}(K)}
\newcommand{\irange}[2]{\llbracket #1,\, #2 \rrbracket}
\newcommand{\Ppkk}[4]{\mathbb{P}_{#1,\,#2}\times \mathbb{P}_{#3,\,#4}}

\newcommand{\SInt}[3]{\int\limits_{#1}^{#2}#3}
\newcommand{\SIntL}[3]{\int_{#1}^{#2}#3}

\usetikzlibrary{calc}

\AtBeginDocument{%
	\setlength\abovedisplayskip{0.4em}
	\setlength\belowdisplayskip{0.4em}}

\newtheorem{property}{Property}
\newtheorem*{possibilities}{Available types of \Dofss}
\newtheorem{conditions}{Requirement}
\newtheorem{conditionsSet}{Admissibility conditions}
 
\newtheorem{example}{Example}
\newtheorem{remark}{Remark}

\newenvironment{Note}{\noindent\textit{\textbf{Note.}}}{\hfill$\blacktriangle$}

\newenvironment{edit}{%
	\leavevmode\ignorespaces%
}{%
}%

\newenvironment{edittt}{%
	\leavevmode\ignorespaces%
}{%
}%

\newenvironment{editt}{%
	\leavevmode\ignorespaces%
}{%
}%

\graphicspath{{Figures/}}
\newdimen\Rr
\usetikzlibrary{calc}
\usetikzlibrary{shapes.misc}
\usetikzlibrary{patterns}
\usetikzlibrary{decorations.pathreplacing}
\usetikzlibrary{matrix,decorations.pathreplacing}
\usetikzlibrary{arrows}
\usetikzlibrary{shapes,positioning}
\tikzset{every picture/.style={remember picture}}
\tikzset{cross/.style={cross out, draw=black, fill=none, minimum size=2*(#1-\pgflinewidth), inner sep=0pt, outer sep=0pt}, cross/.default={2pt}}
\def\centerarc[#1](#2)(#3:#4:#5)
{ \draw[#1] ($(#2)+({#5*cos(#3)},{#5*sin(#3)})$) arc (#3:#4:#5); }

\makeatletter
\pgfdeclarepatternformonly[\LineSpace]{my north east lines}{\pgfqpoint{-1pt}{-1pt}}{\pgfpoint{-\LineSpace+0.2pt}{\LineSpace+0.2pt}}{\pgfqpoint{\LineSpace}{\LineSpace}}%
{
	\pgfsetcolor{\tikz@pattern@color}
	\pgfsetlinewidth{0.4pt}
	\pgfpathmoveto{\pgfqpoint{0pt}{0pt}}
	\pgfpathlineto{\pgfqpoint{-\LineSpace + 0.1pt}{\LineSpace + 0.1pt}}
	\pgfusepath{stroke}
}
\newdimen\LineSpace
\tikzset{
	line space/.code={\LineSpace=#1},
	line space=3pt
}
\makeatother

\newcommand{\SmallMatrix}[2]{
	\usetikzlibrary{matrix,positioning,decorations.pathreplacing}
	\resizebox{!}{\textwidth}{
		\hspace{-15em}\vspace{-30em}
		\begin{tikzpicture}
		[
		style1/.style={
			matrix of math nodes,
			every node/.append style={text width=#1,align=center,minimum height=5ex},
			nodes in empty cells,
			left delimiter=(,
			right delimiter=),
		},
		style2/.style={
			matrix of math nodes,
			every node/.append style={text width=1.4,align=center,minimum height=5ex},
			nodes in empty cells,
			left delimiter=\lbrace,
			right delimiter=\rbrace,
		}
		]
		\matrix (b) [matrix of math nodes,left delimiter=(,right delimiter={)}] {
			\sigma_{I_1}\\
			\phantom{\int}\\
			\sigma_{I_P}\\
		};
		\draw[dashed](b-1-1.south) -- (b-3-1.north);
		\begin{scope}[shift={(8.5,0)}]
		\matrix[style1=1.4cm,ampersand replacement=\&] (1mat)
		{
			\int e_iu_1\cdot p_{1,\,i} \& \& \int e_iu_A\cdot p_{1,\,i} 
			\&\&
			\int x\tilde{u_1}\cdot {p}_{1} \& \& \int x\tilde{u_B}\cdot {p}_{1} \\
			\& \& \& \&
			\& \&
			\\
			\int e_iu_1\cdot p_{P,\,i} \& \& \int e_iu_A\cdot p_{P,\,i}
			\& \&
			\int x\tilde{u_1}\cdot {p}_{P} \& \& \int x\tilde{u_B}\cdot {p}_{P}\\
		};
		\draw[dashed](1mat-1-1.south) -- (1mat-3-1.north);
		\draw[dashed](1mat-1-3.south) -- (1mat-3-3.north);
		\draw[dashed](1mat-1-5.south) -- (1mat-3-5.north);
		\draw[dashed](1mat-1-7.south) -- (1mat-3-7.north);
		
		\draw[dashed](1mat-1-2.center) -- (1mat-1-2.east);
		\draw[dashed](1mat-3-2.center) -- (1mat-3-2.east);
		\draw[dashed](1mat-1-6.center) -- (1mat-1-6.east);
		\draw[dashed](1mat-3-6.center) -- (1mat-3-6.east);
		
		\draw[decoration={brace,raise=5pt},decorate,thick] (1mat-1-1.north west) -- node[above, yshift=0.5em] {$\substack{\text{Block repeats as many time as coordinates;} \\ i\in\irange{1}{d} \\}$} (1mat-1-5.north west);
		\draw[decoration={brace,raise=3pt},decorate,transform canvas={yshift=0.2em},thick] (1mat-1-5.north west) -- node[above] {$\substack{\text{Single block}\\ \phantom{p} }$} (1mat-1-7.north east);
		
		\draw[dashed](1mat-1-4.north east) -- (1mat-3-4.south east);
		\end{scope}
		
		\begin{scope}[shift={(16.2,0)}]
		\matrix (u) [matrix of math nodes,left delimiter=(,right delimiter={)}] {
			a_{i,\,1}\\
			\phantom{\int} \\
			a_{i,\,A}\\
			\\
			b_1\\
			\phantom{\int} \\
			b_P\\
		};
		\draw[dashed](u-3-1.south west) -- (u-3-1.south east);
		\draw[dashed](u-1-1.south) -- (u-3-1.north);
		\draw[dashed](u-6-1.north) -- (u-7-1.north);
		\draw[decoration={brace,raise=15pt},decorate,thick] (u-1-1.east) -- node[right=15pt] {$\substack{\text{Repeats }\\ d \text{ times} }$} (u-3-1.east);
		
		\end{scope}
		
		\draw (1.6,0) node  [align=left] {$ = $};
		
		\end{tikzpicture}
	}
}

\title{General polytopal $\bm{H(\div)}$ -- conformal finite elements and their discretisation spaces}
\author{R\'{e}mi Abgrall, \'{E}lise Le M\'{e}l\'{e}do, Philipp \"{O}ffner}
\date{\today}

\begin{document}

	\maketitle
	
	\begin{abstract}
		We present a class of discretisation spaces and $H(\div)$ -- conformal elements that can be built on any polytope. Bridging the flexibility of the Virtual Element spaces towards the element's shape with the divergence properties of the Raviart -- Thomas elements on the boundaries, the designed frameworks offer a wide range of $H(\div)$ -- conformal discretisations. As those elements are set up through \Dofss, their definitions are easily amenable to the properties the approximated quantities are wished to fulfil. Furthermore, we show that one straightforward restriction of this general setting share its properties with the classical Raviart -- Thomas elements at each interface, for any order and any polytopal shape. Then, to close the introduction of those new elements by an example, we investigate the shape of the basis functions corresponding to particular elements in the two dimensional case.
	\end{abstract}

	\section{Introduction}
\begin{edittt}
There has been recently a lot of activity in the design of methods able to deal with polygonal meshes. In the case of  elliptic PDEs, one can mention  \cite{zbMATH06666913,VEIGA2013,zbMATH06966728,zbMATH06823756,zbMATH06596741}. One may also mention \cite{Talischi2013,zbMATH07078646,zbMATH05897522} for fluid problems, both for compressible and incompressible flow. These are  very partial lists only. In the case of hyperbolic problems, the flexibility of the discontinuous Galerkin method also enables, {\it a priori}, to deal with polygonal meshes, see  \cite{cockburn2012discontinuous}, one of the very first papers after the seminal works of Reeds and Hill \cite{reeds} and Lesaint-Raviart \cite{lesaint}. There are several reasons for this interest of the scientific community: it allows more flexibility in the geometry description and  facilitate mesh adaptation, and several other reasons such as the status of hanging nodes. 

There are many variants of the discontinuous Galerkin method, and one family of algorithms  that has received a lot of attention recently is the so-called Flux Reconstruction methods \cite{huynh2007flux,vincent2011new,ranocha2016summation}. 

The classical simplicial Flux Reconstruction approach involves a point-wise approximation of the flux in a finite-differece framework, modified by a term living only on the boundary in order to retrieve stability and local conservation. Up to some involved modifications, this method can be adapted to quadrangles and hexaheda. In any of those cases, the Flux Reconstruction technique can be rewritten as a Galerkin method applied to a perturbed flux, perturbation guaranteeing the local conservation and stability. The perturbation term is determined by a lifting technique involving a Raviart-Thomas \cite{RTOriginal} polynomial.

Going further and using the Residual Distribution framework, we have shown in  \cite{abgrall:hal-01820176} how to construct schemes analogous to Flux Reconstruction for arbitrary polytopes, convex or non-convex. By using an entropy correction term, this method benefits from a non-linear entropy stability property. 

In this paper, our primary motivation is to construct a Raviart-Thomas like approximation, so that we can reinterpret the correction term we introduce in \cite{abgrall:hal-01820176} exactly as it is done for simplices in \cite{vincent2011new}. We are also interested in providing a new $H(\div)$-conformal discretisation spaces on arbitrary polytopes which we believe has is interest by itself.

\end{edittt}

The theory of $H(\div)$-conformal element has already been studied by  Raviart, Thomas \cite{RTOriginal} and later generalized by N\'{e}d\'{e}lec \cite{Nedelec1980} or Brezzi, Douglas and Marini in the context of mixed finite element method \cite{Brezzi1985}.
More recently, a mixed Petrov-Galerkin scheme using Raviart-Thomas elements has also been investigated in \cite{Dubois2017}. However, up to the authors knowledge, those elements are limited to simplicial and quadrangular shapes.

Several attempts to use general polygons have been made \cite{Gillette2014,Sukumar2006,chen2017minimal},
but they usually make use of generalized barycentric coordinates and are delicate to handle in distorted non-convex elements.
A first polygonal $H(\div)$-conformal element has been proposed in \cite{dipietro2014} using gradient reconstruction and pyramidal sub-meshes tessellation. However, their construction requires some shape regularity within the mesh and the parallel with Raviart-Thomas spaces is limited to the lowest order space. \EEdit{Another approach using stabilisation techniques has been later investigated in \cite{10.1093/imanum/drw003}, where if more flexibility on the element's shape has been achieved the parallel with the Raviart-Thomas elements is still limited to lowest-order simplicial shapes.}

Some other approaches as Virtual Elements Method \cite{VEIGA2013} introduced approximation spaces based on Poisson's solutions. Although more flexible towards the element's shape, those are scalar and not $H(\div)$-conformal. A first promising shape-flexible $H(\div)$-conformal discretisation has been recently proposed in \cite{da2016h}, where the conformity property is enforced directly at the level of the boundaries normal components. 
\Edit{
As it therefore leaves on each boundary face only a scalar representation of the normal component, the boundary setting may appear quite restrictive, especially in applications for which a direct component-wise characterisation or tangential information is suitable.
}
Therefore, its usage may be delicate when one targets a discretisation enhancing the representation of boundary's quantities. We rather focus on creating spaces enhancing a boundary characterisation of the discretised quantity itself, that links to the setting presented in \cite{VEIGA2013}. The quantity of information available on an element is maximized at the boundary, while the $H(\mathrm{div})$-conformity is preserved.

In this paper, we propose a construction enhancing a boundary characterisation, that inherits the interface properties of the Raviart-Thomas elements and benefits from the
shape flexibility of the Virtual Element discretisation. Moreover, rather than defining a basis on which the correction functions can be decomposed, this new setting offers a
new element class that can be used as such in the construction of further numerical schemes. Note also that in \cite{abgrall:hal-01820176}, we do not need the explicit construction of this new elements, we only need the knowledge of degrees of freedom. In that sense, this construction is also in the spirit of VEM.

After briefly recalling the key ideas of the Raviart-Thomas elements, we introduce our new class of discretisation spaces. In a second time, we detail a possible definition of
$H(\div)$-conformal element to finally test them through the behaviour of their corresponding basis functions in the numerical results. For the sake of readability, proofs and further element examples will be given in the appendix. For interested readers, more details on the construction can be found in the extended technical report \cite{abgrall2019class} and an application of those elements can be found in \cite{abgrall:hal-01820176}.

	\section*{Notations}Throughout the paper, our notation will be the following:

\noindent\textit{Geometrical notations:}
\begin{itemize}
\item $x$: spatial variable $x = (x_1,\,\cdots,\,x_d)^T\in\R^d$.
\item $K$: polytopal shape with interior $\mathring{K}$ and  boundary $\partial K$.
\item $\mathfrak{n}$: number of faces of the polytope $K$.
\item $f$: generic face $f$ of the boundary $\partial K$.
\item $n$: generic normal to a generic face $f$.
\item $f_j$: $j^{\mathrm{th}}$ face of the boundary $\partial K$, where $j\in\irange{1}{n}$.
\item $\partial_j  K$: {hyper-face} of $f\in\partial K$ for which the variable $x_j$ is fixed.
\end{itemize}
\textit{Monomials and Polynomial spaces}
\begin{itemize}
\item $\alpha$: multi-index $\alpha = (\alpha_1,\,\cdots,\,\alpha_d)$ defining the monomial $x^\alpha = x_1^{\alpha_1}\cdots x_d^{\alpha_d}$
\item $\Qk$: space of polynomials of degree $\max\limits_i(\alpha_i)$ is at most $k$, of dimension $(k+1)^{d}$
\item $\Pk$: space of polynomials of degree $\sum_{i=1}^{d}\alpha_i\leq k$, of dimension {\scriptsize$\begin{pmatrix}k+d\\k \end{pmatrix}$}
\item $\Qkk$: space of polynomials of degree $\max\limits_i(\alpha_i)=k$
\item $\Pkk$: space of polynomials of degree $\sum_{i=1}^{d}\alpha_i=k$
\item $\mathbb{P}_{k_1,\,k_2,\,\cdots,\,k_d}$: space of polynomials with degree in  $x_i$ $\leq k_i$
\item $\mathbb{P}_{k_1,\,\cdots,[k_i],\,\cdots,\,k_d}$: space of polynomials with degree in  $x_i$ $ k_i$
\item $\mathbb{Q}_{-1} = \mathbb{P}_{-1} = \{0\}$
\end{itemize}
For example, when considering an element in $2D$, a face is one-dimensional and can be parametrised as $x_2 = ax_1+b$ for some constant $a$ and $b$. The space $\partial_1 K$ is then the hyper-face of the line for which the variable $x_1$ is fixed, and $\mathbb{Q}_k(\partial_1 K)$ is reduced to the space of constants. In 3D, $\partial_1 K$ would reduce to the line driven by $x_2$ and  $\partial_2 K$ to the line driven by $x_1$. There, $\mathbb{Q}_k(\partial_1 K)$ would be the space of polynomials of degree $k$ whose monomials are only involving the terms $x_2$.

\noindent\textit{Functional spaces}
\begin{itemize}
\item $\mathcal{R}_k(\partial K) =  \{p\in L^2(K),\,p|_{f_i}\in\Pk(f_i) \text{ for every face } f_i\in\partial K\}$
\item $H(\mathrm{div}, K) = \{u\in \left(L^2(K)\right)^d_,\, \div\, u \in L^2(K)\}$
\end{itemize}
\noindent\textit{Operators}
\begin{itemize}
\item $\bigtimes_{i=1}^{d}$: Cartesian product: $\bigtimes_{i=1}^{d}(x_i)=x=(x_1, \ldots, x_d)$
\item $\{\zeta_i\}_i$: set containing the $d$ cyclic permutations of $\{k+1,\,\underbrace{k,\,\dots,\,k}_{d-1\text{ times }}\}$
\item $\mathbbm{1}$: indicator function
\end{itemize}

	\section{Classical Raviart-Thomas elements.}
The spirit of the Raviart-Thomas elements is to work in a vectorial polynomial discretisation subspace of $H(\mathrm{div}, K)$  in which the functions are characterised separately on the boundary and within the elements. Doing so, the enforcement of the $H(\div)$-conformity can be done at the interfaces by a specific choice of moment-based \Dofs acting only on the boundaries. As formalised by N\'{e}delec \cite{Nedelec1980}, its definition on any simplicial reference shape $K$ contained in $\mathbb{R}^d$ reads
\begin{equation}RT_k(K) = (\Pk(K))^d\oplus x\, \Pkk(K).\label{eq:RTspacetri}\end{equation} 
There, any element $p\in\Rtt \subset (\mathbb{P}_{k+1}(K))^d$ writes under the form
\begin{equation}
	p = \begin{pmatrix}p_1 + x_1\,q\\p_2+x_2\,q\\\vdots\\p_d+x_d\,q	\end{pmatrix}\coloneqq  \bigtimes_{i=1}^{d}(p_i+x_i\,q) = \bigtimes_{i=1}^{d}p_i + \bigtimes_{i=1}^{d}x_i\,q
\end{equation} for some $p_i\in\Pk(K),\, i\in\irange{1}{d}$ and $q\in\Pkk(K)$.
Up to some straightforward computations, the dimension of $RT_k(K)$ can be formulated as
\begin{equation}\label{eq:split}
	\dim \Rtt 
	= \dim (\mathbb{P}_{k-1}(K)^d) + (d+1)\dim \Pk(f).
\end{equation}
Therefore, the definition of the element is done by setting \emph{internal} and \emph{normal} moments projecting respectively on the spaces $\mathbb{P}_{k-1}(K)^d$ and $\mathcal{R}_k(\partial K)$.
\begin{definition}[Degrees of freedom]\label{Def:RtTriDofs}
	Any  $q\in\Rtt$ is determined by
	\begin{subequations}
		\label{eq:TriangleMoments}
		\begin{align}
		\label{eq:TriangleMomentsNormal}
		\textit{Normal moments:} \quad q&\longmapsto\IntL{\partial K}{}{q\cdot n\,p_k}{\gamma (x)}, \quad\forall p_k\in\Rk,\\
		\label{eq:TriangleMomentsInternal}
		\textit{Internal moments:} \quad q&\longmapsto\IntL{K}{}{q\cdot p_{k-1}}{x}, \quad\quad\,\forall p_{k-1}\in(\mathbb{P}_{k-1}(K))^d,
		\end{align}
	\end{subequations}
	where $\mathrm{d}\gamma$ represents the Lebesgue measure on the faces.
\end{definition}
The basis functions of $\Rtt$ that are dual to those \Dofs verify
\begin{equation}
\label{spliit}
q\cdot n|_{\partial K}  \in \Pk(\partial K)  \quad \text{ or } \quad q\cdot n|_{\partial K} \equiv 0
\end{equation}
and are classified respectively as \emph{normal} and \emph{internal} basis functions. One can observe the $H(\div)$-conformity property: it  reduces here to the continuity of the normal component across the boundary. Further divergence properties also hold directly by the nature of the approximation space \cite{RTOriginal}, being a subspace of $H(\div,\,K)$.
\begin{property}[Divergence properties]\label{Prop:DivPropRt}
	For any $q\in\Rtt$, it holds:
	\begin{equation}
	\begin{cases}
	\div\, q \in \Pk(K)\\
	q\cdot n|_{\partial K} \in \Rk.
	\end{cases}
	\label{eq:DivTri}
	\end{equation}
\end{property}
 However, as the relation $(d+1)\dim \Pk(f) = \dim \Rk$ is only valid when the number of edges is $d+1$, this definition is very specific to the simplicial case. Therefore, when going to quads, the definition is changed by modifying the meaning of the polynomial degree, \emph{i.e.} using $\Qk$ spaces instead of $\Pk$ spaces.  The $\Rtt$ space then reads
\[\Rtt = \left(\Qk(K)\right)^d + x\,\Qkk(K) = \bigtimes\limits_{i=1}^{d}\mathbb{P}_{\zeta_i({k+1,\,k,\dots,\,k})}\]\vspace{-0.5em}
and benefits from a dimensional split similar to \eqref{eq:split};
\[\dim \Rtt = 2d \dim \Qk(f) + \dim \bigtimes\limits_{i=1}^{d}\mathbb{P}_{\zeta_i(\{k-1,\,k,\cdots,\,k\})}(K).\]
A  definition of \Dofs analogous to \eqref{eq:TriangleMoments} can then be set up. However, this extension is very specific to quads and cannot be adapted to offer a discretisation framework for arbitrary polytopes (see \cite{abgrall2019class} for details).

	\section{A framework for arbitrary polytopes}
\label{Section:Framework}
In order to build a unifying discretisation framework, we have to define spaces $\Hk$ that fulfil the following property:
\begin{conditions}[Requirements on the discretisation space]
\label{Sec:RequirementsEQGlob}
 The space $\Hk$ is a finite dimensional vectorial subspace of $H(\div,\,K)$ whose dimension adapts to both the number of the polygon's faces and the discretisation order.
\end{conditions}
In addition, to be able to endow $H(\div,\,K)$-conformal elements $E_k(K)$ through definitions of \Dofss, we further ask the {requirement \ref{Sec:RequirementsEQLoc}}.
\begin{conditions}[Requirements on the elements]	\label{Sec:RequirementsEQLoc} \qquad
	\begin{enumerate}
		\item For any space $\Hk$, there exists a unisolvent set of \Dofs $\{\sigma\}$ that can be split into internal and normal subsets so that both the number of internal \Dofs and of the normal \Dofs per face do not depend on the shape of $K$.
		\item The number of internal and normal \Dofs both increase strictly monotonically with the discretisation order.
	\end{enumerate}
\end{conditions}
Lastly, to ensure the existence of a split into internal and normal subsets of \Dofs that matches the classification \eqref{spliit} at the level of the dual basis functions, the feasibility of the {requirement \ref{enum:3}} in $\Hk$ is also needed.
\begin{conditions}[Requirement on the basis functions]
	\label{enum:3}
For any polytope $K$, the internal basis functions vanish on every face of the element.
\end{conditions}
One may possibly ask for one further requirement ensuring a parallel with the Raviart-Thomas setting from the lowest order on.
\begin{conditions}[Optional requirement on the basis functions]
	\label{enum:4}
The lowest order element has no internal \Dofss.
\end{conditions}
\subsection{A class of admissible approximation spaces.}
\subsubsection*{Construction of spaces of discretisation.}
In order to  design a subspace of $H(\div,\,K)$ that satisfies \textit{Condition} \ref{Sec:RequirementsEQGlob} we are led to define a space $\Hk$ with the same architecture as the classical Raviart-Thomas space. Thus, we look for spaces in the form $\Hk = (A_k)^d + x\,B_k$ for two given functional sets $A_k$ and $B_k$.
In order to design those two sets, we start by observing that the use of polynomial spaces is excluded by the condition \ref{enum:3}, being required for any number of edges. Therefore, we consider the spaces $A_k$ and $B_k$ based on solutions to Poisson's problems as in the context of the VEM method \cite{VEIGA2013}. There, a way to allow the existence in $\Hk$ of smooth internal basis functions is to use the set of solutions to the boundary problems $\left\{u|_{\partial K}=0,\,\Delta u=p_k\right\}$ for any $p_k$ belonging to $\mathbb{Q}_m(K)$, $m\in \mathbb{N}\cup\{-1\}$.

In addition, as the $H(\div,\,K)$-conformity will be enforced by normal quantities that are tested only on the boundaries, we also consider the set of Poisson's problems $\{u|_{\partial K} = p_k\mathbbm{1}_{f},\,\Delta u = 0\}$ defined from polynomial boundary functions $p_k\in\mathbb{Q}_l(f)$, $l\in \mathbb{N}\cup\{-1\}$ for each face $f$ of $\partial K$.
Thus, seeing the boundary $\partial K$ face-wise, we define the set
\begin{equation}
\label{eq:Hkcal}
\hk = \{u|_{\partial K}\in L^2(\partial K),\, u|_f \in \Qk(f),\,\, \forall f \in\partial K\}
\end{equation}
and build the space $\Hk$, for integers $l_1,\,l_2,\,m_1$ and $m_2$, {as follows}.
\begin{definition}[$\Hk$ space]
\begin{equation}\begin{split}
\label{eq:HkConstruction}
\Hk =& \{ u \in H^1(K),\,u|_{\partial K}\in \mathcal{H}_{l_1}(\partial K),\, \Delta u \in \mathbb{Q}_{m_1}(K) \}^d \\
&{+}\, x \,\{ u \in H^1(K),\,u|_{\partial K}\in \mathcal{H}_{l_2}(\partial K),\, \Delta u \in \mathbb{Q}_{[m_2]}(K)\}.
\end{split}
\end{equation}
\end{definition}
The choice of $l_1,\,l_2,\,m_1$ and $m_2$ is related to $k$ and will be discussed below.

\begin{remark}\qquad
	\begin{itemize}
	 \item The two subspaces in the definitions of $\mathbb{H}_k(K)$ are in direct sum whenever $l_1\leq0$.
	 \item The presented space is based on polynomial spaces $\mathbb{Q}_k(K)$ rather than $\mathbb{P}_k(K)$ for the sake of consistency with the definition of the Raviart-Thomas space built on quads. This is also a more natural choice when considering mappings to elements of reference, as the monomials involved in the transformations a more coherent with a $\mathbb{Q}_k(K)$ based discretisation {(especially for the lower order discretisation where the monomial $xy$ is not part of $\mathbb{P}_1(K)$ but already belongs to $\mathbb{Q}_1(K)$)}.
	\end{itemize}

\end{remark}

\subsubsection*{Properties of $\Hk$ spaces.}\label{Sec:HkProperties}
The space $\Hk$ is constructed from four independent blocks whose definitions are driven by the independent coefficients $l_1,\,l_2,\,m_1$ and $m_2$.  The couple $(m_1,\,m_2)$ drives the discretisation quality exclusively within the cell while $(l_1,\,l_2)$ takes care only of the boundary. Thus, the separation between internal and normal basis functions is natural. Furthermore, the {Property \ref{Prop:Hdiv}} holds, emphasising that the $H(\div,\,K)$-conformity is ensured by the definition of $\hk$, while the inner smoothness is provided through the Laplacian.
\begin{proposition}
	\label{Prop:Hdiv}
		For any function $q$ belonging to any space $\Hk$, it holds:
\begin{equation}\label{Prop:ConformityHk}
		q\,\cdot n |_{\partial K} \in \mathcal{H}_{\max\{l_1,\,l_2\} }(\partial K)
\quad \text{and} \quad
		\div\, q \in L^2(K).
		\end{equation}
\end{proposition}
It comes the following inclusion allowing $H(\div)$-conformity.
\begin{corollary}\label{Prop:NgonsConformity}
		For any couples  $(l_1,\,l_2)$ and $(m_1,\,m_2)$,
		\[
		\Hk \subset H(\div,\,K).
		\]
\end{corollary}
When selecting $l_1\leq0$, those spaces are of dimension 
\begin{equation}
\label{Prop:DimNgons}
\dim \Hk \,=\, {\mathfrak{n}} \big (d(l_1+1)^{d-1} + (l_2+1)^{d-1}\big)\,\, + \,\,\big(d(m_1 + 1)^d + (m_2 + 1)^d - m_2^d\big),
\end{equation}
making their structure {a-priori} suitable to be used as discretisation spaces endowing $H(\div,\,K)$-conformal elements.
\begin{example}
	When $K$ is a two-dimensional simplex and when $l_1$, $l_2$ are chosen as $(l_1,\,l_2)=(-1,\,k)$, the discretisation quality of the normal component matches the one of the Raviart--Thomas setting.
\end{example}
\subsubsection*{Admissibility of the spaces for building $H(\div)$-conformal elements.}
In order to define elements in the spirit of Raviart-Thomas, we need to set $(d(l_1+1)^{d-1} + (l_2+1)^{d-1})$ normal \Dofs per face and $d(m_1 + 1)^d + (m_2 + 1)^d - m_2^d$ internal \Dofss. While this splitting does not impact the set of admissible coefficients $(m_1,\, m_2)$, it reduces the range of coefficients $(l_1,\,l_2)$ that can be used. Indeed, the space $\Hk$ is constructed from four independent blocks providing two distinct discretisations: on the boundary and within the element. Thus, when testing a function of $\Hk$ through normal \Dofss, one can only retrieve the polynomial obtained from the two boundary conditions defining the sets $A_k$ and $B_k$.
On each face, this polynomial is of the form $p = p_{k,\,A} + p_{k,\,B}$, where the function $p_{k,\,A}\in(A_k)^d|_{f}$ reads
\begin{equation}
\label{eq11}
p_{k,\,A} = \bigtimes\limits_{j = 1}^{d} \Big(\sum\limits_{|\alpha_{i}|\leq l_1} a_{ij}x^{\alpha_{i}}\Big)
\end{equation}
for a given set of multi-index $\{\alpha_{i}\}_{i}$ and coefficients $\{a_{i,\,j}\}_{i,\,j}$  depending on the coordinates $x_j$. The function $p_{k,\,B}\in x\,B_k|_f$ reads however
\begin{equation}
\label{eq12}
p_{k,\,B} = \bigtimes\limits_{j = 1}^{d} \Big(x_j\sum\limits_{|\beta_i|\leq l_2} b_i x^{\beta_i}\Big)
\end{equation}
for a given set of multi-indices $\{\beta_i\}_i$ and coefficients $\{b_i\}_i$ independent of the coordinates $x_j$. Therefore, denoting by $\{\xi_j\}_{j\in\irange{1}{d}}$ the coordinates permutation that allows to shift the lowest orders terms of $x\,B_k|_f$ to $(A_k)^d|_f$, $p\in \Hk|_f$ can be written as follows.
\begin{align}
\label{eq:FormOk}
\text{If}\,\, l_2\geq l_1,\,\,\quad\,\,\,\,\quad p &= \bigtimes\limits_{j = 1}^{d} \Big(\sum\limits_{|\alpha_{i}|\leq l_1} a_{ij} x^{\alpha_{i}} \Big) + \bigtimes\limits_{j = 1}^{d} \Big(x_j\sum\limits_{|\beta_i|\leq l_2} b_{i} x^{\beta_{i}} \Big)\nonumber\\
&= \bigtimes\limits_{j = 1}^{d} \Big(\sum\limits_{\substack{|\alpha_{i}|\leq l_1 \\ |\alpha_{i}|\neq 0}} (a_{ij} + b_{\xi_j(i)})x^{\alpha_{i}} + x_j\hspace{-1em} \sum\limits_{l_1\leq |\beta_{i}|\leq l_2}\hspace{-1em} b_{i} x^{\beta_{i}}  \Big) + \bigtimes\limits_{j = 1}^{d} a_{0j}x_j^0.\\
\text{If}\,\,l_1\geq l_2+1,\quad
 p& = \bigtimes\limits_{j = 1}^{d} \Big(\sum\limits_{\substack{|\alpha_{i}|\leq l_1\\|\alpha_{i}\neq 0|}} (a_{ij} + b_{\xi_j(i)}) x^{\alpha_{i}} + a_{0j}x_j^0\Big).\label{eq:FormOk2}
\end{align}
The structure of those relations implies that the terms $a_{ij}$ and $b_{\xi_j(i)}$ are  combined into a single coefficient and cannot be specified individually from further normal \Dofss. Indeed, the remaining freedom can only be {seen inside the polytope, as a consequence of the boundary conditions} on the Poisson's solutions in either $A_k$ or $x\, B_k$. To prevent any over-determination by the normal \Dofs in $\Hk|_{f}$, we therefore have to make sure that the dimension of the boundary part \eqref{eq:FormOk}-\eqref{eq:FormOk2} of any function living in $\Hk$ is larger than the number of wished normal \Dofs per face. By reading out the structure of \eqref{eq:FormOk}-\eqref{eq:FormOk2} it comes
\begin{equation}\label{eqrestrict}
\dim \Hk|_{f} = \begin{cases}
{(l_2+1)^{d-1}}&{\text{ if }\,  l_1= -1,}\\
d(l_1+1)^{d-1} + (l_2+1)^{d-1} - l_1^{d-1} &\text{ if }\,  l_2\geq l_1,\\
d(l_1+1)^{d-1} &\text{ otherwise}.
\end{cases}
\end{equation}
We thus restrict the admissible couples $(l_1,\,l_2)$ to those verifying the admissibility condition \ref{prop:Cond1}, preventing any over-determination.
\begin{conditionsSet}[Necessary condition for using conformal elements]
	\label{prop:Cond1}
If $\dim \mathcal{N}$ is the number of normal moments per face that we wish, and $\dim \Hk_{|f}$ is the number of coefficient we can tune for the face $f$, we should have: 
$$\dim \mathcal{N}\leq \dim \Hk_{| f}.$$
	In the case $l_2\geq l_1$, it reduces to:
	\[
	d(l_1+1)^{d-1} + (l_2+1)^{d-1} \leq  d(l_1+1)^{d-1} + (l_2+1)^{d-1} - l_1^{d-1} \quad\quad\quad (\Leftrightarrow\quad	 l_1^{d-1} \leq 0)
	\]
	while otherwise it comes
	\[
	d(l_1+1)^{d-1} + (l_2+1)^{d-1} \leq d(l_1+1)^{d-1}\quad\quad\quad\quad\quad\quad\quad\quad\quad\quad\quad\quad (\Leftrightarrow \quad l_2 = -1).
	\]
\end{conditionsSet}
Regarding the internal characterisation, any couple of coefficients $(m_1,\,m_2)$ is allowed. 

\subsubsection*{Definition of series of spaces.}
While fulfilling the above conditions, one can set a specific discretisation framework within which the spaces share a predefined structure. By example, defining the four coefficients $l_1,\,l_2,\,m_1$ and $m_2$  through affine relations of the type $l = ak +b$ for some index $k\in\mathbb{N}$, the range of discretisation qualities achievable within the framework is predetermined by a refinement sequence in each block, and the order of each space can be simply defined as the index $k$ generating each of the four coefficients. A typical working example is obtained by defining $m_1=m_2=k-1,\,l_1=0$ and $l_2=k$, leading to a series of discretisation spaces of order $k$. This case is specifically detailed in the section \ref{Sec:2Dcase}.

\subsection{Definition of admissible elements.}
\label{Sec:admissibleHk}
Under the admissibility conditions \ref{prop:Cond1}, the spaces $\Hk$ allow the construction of $H(\div)$-conformal elements through the definition of normal \Dofs enforcing the conformity and internal ones preserving it. We propose here a possible construction of such sets.

\subsubsection*{Definition of admissible normal \Dofss.}
\label{AdmissibleDofs}
The role of the normal \Dofs is to determine vectorial polynomials on the boundaries and to enforce the $H(\div)$-conformity of the element. We define them as the normal component of the tested quantities projected against polynomials of $\Hk_{| f}$. We focus on the following possibilities:
\begin{possibilities}
	\label{enum:AvailDofsTypes}
	For any $q\in\Hk$, we define:
	\begin{enumerate}
		\item The face integral of coordinate-wise components tested against polynomials:
		\begin{subequations}
			\label{fullset}{\small
		\begin{equation}
				q \mapsto \int_{f}{q_i\,n_{ix}\,p\,}\dd{\gamma(x)}, \quad\, \forall\,p\in\mathbb{Q}_{\max\{l_1,\,l_2\}}(f),\quad\quad\quad\quad\quad\quad\,\,\,\quad\quad\quad\quad\quad\label{1a}
		\end{equation}\par}
		\item The face integral of a function in $\Hk$ projected onto the face normal, and tested against polynomials:{\small
		\begin{align}
		q &\mapsto \int_{f}{q\cdot n\,p\,}\dd{\gamma(x)}, \,\, \forall\,p\in\mathbb{Q}_{\max\{l_1,\,l_2\}}(f),\label{1b}\\
		q &\mapsto \int_{f}{q\cdot n\,p\,}\dd{\gamma(x)}, \,\, \forall\,p\in \{x_i\,p_i, \,p_i\in\mathbb{Q}_{\zeta_i([l_2],\,l_2,\,\cdots,\,l_2)}(f),\,i\in
		\irange{1}{d-1}\}\label{1c}
		\end{align}\par}
		\item The pointwise values of the discretised quantity tested against the face's normal:
		\begin{equation}
		q \mapsto q(x_{im})\cdot n_i,\, \text{ for sampling points } \{x_{im}\}_m \text{ on the face } f_i.\quad\quad\quad	 \label{1d}
		\end{equation}
	\end{subequations}
	\end{enumerate}
\end{possibilities}
Defining the normal \Dofs then reduces to choosing $d(l_1+1)^{d-1} + (l_2+1)^{d-1}$ of them among the possibilities \eqref{fullset} so that their set is unisolvent for $\Hk|_f$. To ensure this, preventing any under-determination is sufficient. Therefore, we need to avoid the selection of projectors that are linearly dependent, and pay attention to determining both global and coordinate-wise behaviours of any vector polynomial $q\in\Hk_{| f}$.
\begin{example}
	\label{ex:Struct}
	In two dimensions and for $l_1=l_2=0$, any $p\in\Hk_{| f}$ reads
	\[\small{q = \begin{pmatrix} A \\ B \end{pmatrix} + C \begin{pmatrix}x\\y\end{pmatrix}}\]
	for some constants $A$, $B$ and $C$. 
	The characterisation of $q$ can be done by selecting two component-wise moments involving  $A\,n_{ix}$ or $B\,n_{iy}$ tested against the constant polynomial $p = 1$ and one global moment that tests $q\cdot n = C\,(n_{ix}+n_{iy})+A\,n_{ix}+B\,n_{iy}$ against the polynomial $p = x$. One could also choose two global moments and one coordinate-wise.
\end{example}
In practice, the selection of \Dofs reduces to choosing the polynomials $p$ on which the function $q$ will be tested coordinate-wise. The other polynomials $p$ play the role of test functions for the global normal component $q\cdot n$. The unisolvence of the set is then ensured by the following admissibility conditions.
\begin{conditionsSet}
$ $
	\label{Ass:Dofs}
	\begin{enumerate}
	\item \label{Ass:Dofs0}
          The projection polynomials $p$, and all the polynomials $\sigma\colon q \mapsto \sigma(q(x_{im}))$ that define the point values must be  linearly independent.
		\item \label{restriction2}
		When using a coordinate-wise degree of freedom of the type \eqref{1a}, polygonal shapes $K$ containing a face parallel to any axis are not allowed. The term $n_{ix}$ or $n_{iy}$ would indeed {always} vanish for some $i\in\irange{1}{n}$, thus not {describing} any function of $\Hk_{| f}$.
	\end{enumerate}
\end{conditionsSet}

\begin{editt}
	\begin{Note}
		The second point of the admissibility conditions may seem unreasonable as it may prevent the use of some shapes for specific orientations. However, it is always possible to easily modify the incriminated moments element-wise or to select other moments that make the element robust with respect to rotation while still yielding $H(\div)$-- conformity. See \cite{abgrall2019class} for more details.
	\end{Note}
\end{editt}
\vspace{0.3em}

To help the construction of an element on $K$ through the selection of \Dofs among those fulfilling the admissibility conditions \ref{Ass:Dofs}, we recall that the chosen set of \Dofs imposes the shape of the dual basis functions. We can therefore select the \Dofs depending on the wished properties of the basis functions.

More crucially, the selection of global and/or coordinate-wise normal \Dofs leads to the reclassification of some basis functions as internal ones. Indeed, as the face-wise normal component of any function $q$ in $\Hk_{| f}$ is only of degree $\max\{l_1,\,l_2\}$, the term $q\cdot n|_{f}$ requires only  $(\max\{l_1,\,l_2\}+1)^{d-1}$ basis functions to be decomposed on. Therefore, up to $d\,(l_1+1)^{d-1}$ basis functions may see their global normal component vanishing on every face. Their coordinate-wise components will however not vanish, as they take care of the coordinate-wise behaviours that cannot be determined solely through the expression of $q\cdot n|_{f}$.

\begin{remark}
	\label{reducing}
	Typically, the more global \Dofs are designed, the more the representation of $p\cdot n$ is completed globally. As a consequence, more basis functions have a vanishing normal component as they are forced to take care only of coordinate-wise behaviours, forcing them to be reclassified into internal basis function. The reverse scenario may also be considered.
\end{remark}

To avoid this reclassification and allow a parallel with the Raviart-Thomas elements from the lowest order space, we ask the {requirement \ref{enum:4}}. This will be discussed in the section \ref{Reduced}.

\subsubsection*{An example of a possible definition of normal \Dofss.}
As an example, we detail one selection of normal \Dofs in the case $l_2\geq l_1$ where every function in $\Hk|_{f}$ is of the form \eqref{eq:FormOk}. For interested readers, other possibilities are presented in the technical report \cite{abgrall2019class}.

Here, we select moments from the set \eqref{fullset} so that the elements of $\Hk_{| f}$ are determined as much as possible by testing only their normal component. The remaining freedom is characterised by few coordinate-wise moments. We consider:
\begin{subequations}
	\label{eq:Doffs}
	\begin{align}
	\sigma\colon q  &\mapsto \IntL{f_j}{}{q_i \,n_{x_i} \,x_i^{l_1+1}}{\gamma(x)},\,&&{\text{ for all } j\in\irange{1}{\mathfrak{n}} \text{ and all } i\in\irange{1}{d},}\label{eq:Doffs3}\\
	\sigma\colon q &\mapsto \IntL{f_j}{}{
		q \cdot \,n\, p_k}{\gamma(x)},  && {\text{ for all } j\in\irange{1}{\mathfrak{n}} \text{ and all } p_k\in\mathbb{Q}_{l_2}(f_j)\setminus\mathbb{Q}_{l_1}(f_j),} \label{eq:Doffs2}\\
	\sigma\colon q &\mapsto \IntL{f_j}{}{q\cdot n \,x_j\, x_j^{l_2} \tilde{x}}{\gamma(x)},\,&& {\text{ for all } j\in\irange{1}{\mathfrak{n}} \text{ and any }  \tilde{x}\in \mathbb{Q}_{l_2}(\partial_j K),}  \label{eq:Dofs111}
	\end{align}
\end{subequations}
where $\tilde{x}\in\mathbb{Q}_{l_2}(\R^{d-2})$ is not involving the variable $x_j$ so that the moment \eqref{eq:Dofs111} has for integrands the second terms of the right hand side of \eqref{eq:FormOk} when $|\beta_i|=l_2$. Note that the set \eqref{eq:Doffs} is of dimension {\small $d(l_1+1)^{d-1}  + (l_2+1)^{d-1} - (l_1+1)^{d-1} + (d-1)(l_2+1)^{d-2}$} though we require {\small $d(l_1+1)^{d-1}  + (l_2+1)^{d-1}$} moments. Thus, this configuration can only be used when $l_1$ and $l_2$ verify the feasibility condition:
\begin{equation}
\label{ineq2}
(l_1+1)^{d-1}\leq (d-1)(l_2+1)^{d-2}_,
\end{equation}
which is a reduction of the {admissibility conditions \ref{prop:Cond1}}.
In two dimensions the above relation reduces to an equality, and all the \Dofs presented in \eqref{eq:Doffs} are considered. In higher dimensions, a further selection from the set \eqref{eq:Doffs} is required. There, we consider the sets \eqref{eq:Doffs3}-\eqref{eq:Doffs2} fully and select any $(l_1+1)^{d-1}$ moments from \eqref{eq:Dofs111}.

\begin{definition}\label{def:def3}
	Any choice of {\small{$(l_1+1)^{d-1}$}} moments among \eqref{eq:Dofs111} is denoted as the \enquote{configuration \IIa}. Associated with any admissible internal \Dofss, its unisolvence is given by the lemma \ref{LemmaImportant}.
\end{definition}

Up to the additional coordinate-wise moments, the configuration \IIa{} is close to the Raviart-Thomas setting. However, the scaling of the dual basis functions does not match the one of the Raviart-Thomas basis. In order to obtain a similar scaling, one should rather scale the above \Dofs with respect to each edge's length and orientation, or consider in place of the moments \eqref{eq:Dofs111} the point-wise values
\begin{equation}
q\mapsto q(x_{im})\cdot n,  \label{eq:Dofs12}
\end{equation}
where $ i\in\irange{1}{d}$, $m\in\irange{1}{(1+l_1)^{d-1}}$ and $x_{im} $ is any sampling point on the face $f_i$.
\begin{definition}\label{def:def4}
	{Any selection of {\small{$(l_1+1)^{d-1}$}} \Dofs among the sets \eqref{eq:Doffs3},\eqref{eq:Doffs2} and \eqref{eq:Dofs12} is labelled as the \enquote{configuration \IIb}. Associated with any admissible internal \Dofss, its unisolvence is given by the {lemma \ref{LemmaImportant}}}.
\end{definition}

\begin{remark}
As we assume in this example that $l_2\geq l_1$, the choice of $l_1$ is restricted to either $l_1=0$ or $l_1=-1$. Thus, in the \emph{definitions} \ref{def:def3} and \ref{def:def4}, one would only need to select respectively one or none \Dofs from the set \eqref{eq:Dofs111} or\eqref{eq:Dofs12}.
\end{remark}

\subsubsection*{Definition of admissible internal \Dofss.}
\label{internal}
In order to define admissible internal \Dofss, we have to make sure that the corresponding internal basis functions vanish on every face. We therefore stick to the idea of Raviart-Thomas and define moment based degrees of freedom that read for any $q \in \Hk$
\begin{equation}
\label{eq:DofsIntHk}
\sigma(q) \mapsto \IntL{K}{}{q \, \cdot \, p_k}{x}, \quad \text{ for all } p_k\in\mathcal{P}(K)
\end{equation}
for some function space $\mathcal{P}(K)$ of dimension $\left((m_1 + 1)^d + (m_2 + 1)^d \right)$.
Considered as a test space, $\mathcal{P}(K)$ may simply gather polynomial functions used in the definition of the Poisson's problems generating $\Hk$. The discretised quantities would then be determined through their polynomial projections. Another choice is to test against the set of Poisson's solutions to the problems $\{\Delta p_k \in \mathcal{P},\, p_k|_{\partial K}=0\}$.

Using one or the other possibility for $\mathcal{P}(K)$, the unisolvence of the set of internal \Dofs in $\Hk|_{\mathring{K}}$ is ensured by the following admissibility conditions (see the proof \ref{LemmaImportant}, Part $3$):
\begin{conditionsSet}\label{Ass:Dofs1}
$ $
	\begin{enumerate}
		\item The polynomials $\{p_l\}_l$ generating $\mathcal{P}(K)$ are linearly independent.
		\item No polynomial $p_l$ is of degree larger than $\max\{m_1,\,m_2+1\}$.
	\end{enumerate}
\end{conditionsSet}

\subsubsection*{Definition of the elements.}
\label{susbsub:further}
Combining the two previous paragraphs with the definition of the space $\Hk$, $H(\div,\,K)$-conformal elements can be set up.
\begin{proposition}
\label{GeneralProof2}
Let $K$ be any polytope satisfying the second item of the admissibility conditions \ref{Ass:Dofs} and $\Hk$ be any admissible space built on it. Let also $\{\sigma_{N}\}$ be any selection of $d(l_1+1)^{d-1} + (l_2+1)^{d-1}$ \Dofs from the set \eqref{fullset} fulfilling first item of the admissibility conditions \ref{Ass:Dofs}, and $\{\sigma_{I} \}$ the set of internal moments built through the expression \eqref{eq:DofsIntHk} for any of the projection sets $\mathcal{P}(K)$ fulfilling the admissibility conditions \ref{Ass:Dofs1}. Then, the set $\{\sigma_{N}\}\cup\{\sigma_{I}\}$ is unisolvent for $\Hk$ and defines a $H(\div,\,K)$-conformal element.
\end{proposition}
This well-possessedness property is an immediate corollary of the following proposition, proven in the appendix \ref{propro}.
\begin{proposition}
  \label{Prop:NgonsUnisolvence3}
  Let $q\in \Hk$, and denote $\sigma_{N}(q)$ the n-tuple of normal \Dofs extracted from the set \eqref{fullset}. If
	\begin{equation*}
	\sigma_{N}(q) = 0\,\, \text{ and }\,\,
	\IntL{K}{}{q \, \cdot \, p_k}{x} = 0 \quad \text{ for all } p_k\in\mathcal{P}
	\end{equation*}
        then $q=0$.
\end{proposition}
At this point, any admissible definition leads to $H(\div,\,K)$-conformal elements.

\subsection{Summary of the construction.}
Let us summarize the spaces construction and the example of normal \Dofs that has been detailed above. To begin with, the class of discretisation spaces reads
\begin{align*}
\Hk =& \{ u \in H^1(K),\,u|_{\partial K}\in \mathcal{H}_{l_1}(\partial K),\, \Delta u \in \mathbb{Q}_{m_1}(K) \}^d \nonumber \\
&\bigoplus x \,\{ u \in H^1(K),\,u|_{\partial K}\in \mathcal{H}_{l_2}(\partial K),\, \Delta u \in \mathbb{Q}_{[m_2]}(K)\},
\end{align*}
with the convention that $\mathbb{Q}_{-1} = \{0\}$ and where the integers $l_1,\,l_2,\,m_1$ and $m_2$ verify
\[
m_1,\,m_2,\, l_2\geq -1\quad\quad \text{ and }\quad\quad -1 \leq \,\, l_1 \leq \,  0.
\]
So defined, it holds:
\begin{itemize}
	\item $\dim \Hk = {\mathfrak{n}}(d(l_1+1)^{d-1}+(l_2+1)^{d-1}) + ((m_1+1)^d+(m_2+1)^d-m_2^2)$
	\item For all $q\in\Hk$, $q\cdot n|_{\partial K} \in \mathcal{H}_{\max\{l_1,\,l_2\}}(\partial K)$ and $\Hk\subset H(\div,\,K)$.
\end{itemize}
Thus, conformal elements can be defined through normal \Dofs enforcing the $H(\div)$-conformity and internal ones preserving it, provided that the polytope $K$ satisfies the two conditions
\begin{itemize}
\item The polytope $K$ has a reasonable aspect ratio, so that the Poisson problem ({required by the definition of the underlying VEM spaces}) is well posed.
\item {No face is parallel to any axis, ensuring the unisolvence of the presented \Dofs (when using component-wise \Dofss).}
\end{itemize}

Note that when selecting component-wise \Dofss, the above condition on the orientation of the face with respect to the axis raises stability issues when dealing with element whose faces are almost parallel to the axis. This issue is easily avoidable by selecting at least a global degree of freedom involving the term $p\cdot n$ in the moment's integrand and changing the testing vector in the coordinate-wise \Dofs to any vector $v\neq n$.
\vspace{1em}

\noindent \textbf{Internal \Dofss.}
It is set
\begin{equation}
\label{eq:DofsIntHk2}
\sigma\colon q \mapsto \IntL{K}{}{q \, \cdot \, p_k}{x}, \quad \forall p_k\in\mathcal{P}
\end{equation}
for any space $\mathcal{P}$ defined either as a polynomial space or as any subspace of Poisson's solutions, having for dimension $(m_1+1)^d + (m_2+1)^d-m_2^d$ and fulfilling the {assumptions \ref{Ass:Dofs1}}.
\newline

\noindent \textbf{Normal \Dofss.}
\begin{table}[h]
\begin{center}\resizebox{\linewidth}{!}{\renewcommand{\arraystretch}{0.5}
	\begin{tabular}{|p{22mm}|p{17mm}|p{20mm}|p{28mm}|}
		\hline
		{\tiny Representation } & {\tiny of low order}  & {\tiny of $(A_k\cap B_k)|_{\partial K}$} & {\tiny of higher orders}\\
		\hline
	\tiny{Available when \phantom{fwpo}}	$  \scriptsize{(d-1)(l_2 + 1)^{d-2}}$\newline $ \tiny{\geq (l_1+1)^{d-1}_.}$\newline $ $ \newline  Select $(l_1+1)^{d-1}$ moments per face from the bold ones.& \tiny{Inherited from\newline the highest order representation} & \vspace{-0.2em} {$ $		  {\tiny{$\SInt{f_j}{}{q_i \,n_{x_i} x_i^{l_1+1}}$}}}\newline$ $ \newline $ $  \newline $ $ \newline $ $
		\newline \tiny{$\forall i\in\irange{1}{d},\,$}\newline \tiny{{$\forall j\in\irange{1}{\mathfrak{n}}$}}
		& {{\tiny$\SInt{f_j}{}{q \cdot \,n\, p_k} $}}
		\newline {{\tiny$\boldsymbol{\SInt{f_j}{}{q\cdot n x_i x_i^{l_2} \tilde{x}}}$}}
		\newline\newline\newline {\tiny {$ \forall p_k\in\mathbb{Q}_{l_2}(f_j)\setminus\mathbb{Q}_{l_1}(f_j),\, $}}
		\newline\tiny {{$\forall\tilde{x}\in \mathbb{Q}_{l_2}(\partial_i K)$}}\\\hline
	\end{tabular}}
\end{center}
\caption{Summary of the used \Dofs for the configurations \IIa.\label{MomentsRecap0}}
\end{table}
\noindent 
Though in two dimensions the above setting is fixed and all the mentioned \Dofs have to be considered, in three dimension the selection of \Dofs among the bold ones is a matter of taste, possibly directed by properties of the discretised quantities that are known \emph{a priori}. Note also that one could project on any other polynomial basis rather than using projections over monomials.

\FloatBarrier

\subsection{Two examples in two dimensions.}
\label{Sec:2Dcase}

We first detail  an example of a discretisation framework contained in the previously presented setting for which a parallel with the
Raviart-Thomas elements can be drawn from the order $k=1$ on. In a second time, we present an example of a reduced framework where a parallel with the Raviart-Thomas is achieved at any order.

\subsubsection{An example of a general setting.}
\label{Sec:2Dcase1}
We consider a series of discretisation spaces by indexing the coefficients $l_1 = 0,\, l_2 = k,\, m_1 = k-1$ and $m_2 = k-1$ for any $k\in\mathbb{N}$, seen here as the space order. The space $\Hk$ is then defined as
\begin{align}
\label{eq:TheSpace2D}
\Hk =& \{ u \in H^1(K),\,u|_{\partial K}\in \mathcal{H}_{0}(\partial K),\, \Delta u \in \mathbb{Q}_{k-1}(K) \}^2\nonumber \\
&\oplus x \,\{ u \in H^1(K),\,u|_{\partial K}\in \mathcal{H}_{k}(\partial K),\, \Delta u \in \mathbb{Q}_{[k-1]}(K)\}.
\end{align}
By a straightforward application of the previous section, it comes
\begin{equation}
	\dim \Hk = \mathfrak{n}(k+3) + 2k(k+1) - \mathbbm{1}_{k>0}.
\end{equation}

\subsubsection*{Example of a two dimensional element.}
The example of selected normal \Dofs defining the elements $E_k = (K,\, \Hk,\,$ $ \{\sigma\})$ presented in the previous section then reduces to the expressions given in the table
\ref{eq:DofsIntHk2D1}.
\begin{table}[h]
	\begin{center}\resizebox{0.8\linewidth}{!}{\renewcommand{\arraystretch}{0.5}
			\begin{tabular}{|p{22mm}|p{22mm}|p{30mm}|p{30mm}|}
				\hline
			 {\scriptsize Representation} & {\scriptsize of low order }  & {\scriptsize of $(A_k\cap B_k)|_{\partial K}$} & {\scriptsize of higher orders}\\
				\hline
			 {\scriptsize Moments }&{\small$\SInt{f_j}{}{q\cdot n}$}\newline$ $\newline$ $\newline
			 	\scriptsize{$\forall j\in\irange{1}{\mathfrak{n}}$}
				& \small{$\SInt{f_j}{}{q_{i} \,n_{x_i} x_i}$}
				\newline $ $\newline $ $\newline
				{\scriptsize{$\forall i\in\irange{1}{2},\,\forall \,j\in\irange{1}{\mathfrak{n}}\newline$}}
				& {\small$\SInt{f_j}{}{q \cdot \,n\, p_k}$} \newline$  $\newline $ $\newline
				{\scriptsize{$\forall j\in\irange{1}{\mathfrak{n}},$}}\newline
					{\scriptsize{$\forall p_k\in\mathbb{Q}_{k}(f_j)\setminus\mathbb{Q}_{0}(f_j)$}}\\\hline
		\end{tabular}}
	\end{center}
	\caption{Definition of the degrees of freedom in the 2D case for the configuration \IIa.\label{eq:DofsIntHk2D1}}
\end{table}
\noindent The internal \Dofs are set as
\begin{equation}
\label{eq:DofsIntHk2D2}
\sigma(q) \mapsto \Int{K}{}{q \, \cdot \, p_k}{x}, \quad \text{ for all } p_k\in\mathcal{P},
\end{equation}
where $\mathcal{P}$ is chosen as the symmetric space
 \[\mathcal{P} = \Ppkk{k}{k-1}{k-1}{k}\setminus \left(\Ppkk{[k]}{[k-1]}{[k-1]}{[k]}\right)\cup \left\{(x,\,y)^T\mapsto \begin{pmatrix}x^ky^{k-1} \\ x^{k-1}y^k\end{pmatrix}\right\}.\]
  Though the internal projection space is less refined than the one set on the edges, this is not bothersome as the impact of the divergence within the cell is less dramatic. Note also that in practice, for defining the projections \emph{\eqref{eq:DofsIntHk2D2}} one can work with any basis of $\mathcal{P}$.

\subsubsection*{Link to another class of elements.}

As pointed out in the introduction, this contribution  can be linked to a discretisation setting presented in \cite{da2016h}, where the considered space reads
\begin{align*}
   \mathcal{V}_{2,\,k}^{face}(K) = \{v\in& H(\mathrm{div},\, K)\cap H(\mathrm{curl},\, K)\quad s.t. \quad v\cdot n|_f \in \mathbb{P}_k(f) \forall f \in \partial K,\\
   & \mathrm{grad} (\mathrm{div} (v)) \in \nabla \mathbb{P}_{k-1}(K),\, \text{and } \mathrm{curl} \,v \in \mathbb{P}_{k-1}(K)\}.
\end{align*}
\begin{edit}
Restricting the setting on the boundary to polynomial functions, that is introducing 
\begin{align*}
\tilde{\mathcal{V}}_{2,\,k}(K) = \{v\in& H(\mathrm{div},\, K)\cap H(\mathrm{curl},\, K) \quad s.t. \quad v\cdot n|_f \in \mathbb{P}_k(f) \, and \, v|_f \in\bigtimes_{i=1}^d\mathbb{Q}_{\zeta_i(k+1,\,k,\,\cdots,\,k)}(f)\quad \forall f \in \partial K,\\
& \mathrm{grad} (\mathrm{div} (v)) \in \nabla \mathbb{P}_{k-1}(K),\, \text{and } \mathrm{curl} \,v \in \mathbb{P}_{k-1}(K)\}.
\end{align*}
\end{edit}
it can in particular be shown that
\begin{equation}\label{eq:connection}
 v\in \Edit{\tilde{\mathcal{V}}}_{2,\,k}\quad  \Rightarrow \quad v\in \mathbb{H}_k(K),
\end{equation}
for the space $\mathbb{H}_k(K)$ is constructed from the coefficients $(l_1,\,l_2) = (0,\,k)$ and $(m_1,\,m_2) = (k,\,-1)$ as
\begin{equation*}\begin{split}
\Hk =& \{ u \in H^1(K),\,u|_{\partial K}\in \mathcal{H}_{0}(\partial K),\, \Delta u \in \mathbb{Q}_{k}(K) \}^d \\
&{+}\, x \,\{ u \in H^1(K),\,u|_{\partial K}\in \mathcal{H}_{k}(\partial K),\, \Delta u \in \mathbb{Q}_{[-1]}(K)\}.\end{split}
\end{equation*}
Indeed, any element $v\in\tilde{\mathcal{V}}_{2,\,k}(K)$ belongs to $H^1(K)$, and can be written as $v = v_1 + v_2$, where $v_1$ lives in $A_k$ and $v_2$ lives in $x\,B_k$ (see the \emph{Appendix} for a sketch of the proof).

The \Dofs selected in the framework of \cite{da2016h} are however different from what we do here. We allow more freedom in the inner characterisation while preserving the desired properties on the boundary. Note also that contrarily to the more general setting presented in this section, the lowest order elements of  \Edit{\cite{da2016h} cannot be natively defined. As the normal component on the boundary belongs at least to  $\mathbb{P}_1(K)$ for any admissible space, the setting of \cite{da2016h} has to be slightly modified (see e.g. \cite{refId0})}.

\subsubsection{An example of a reduced setting.}
\label{Reduced}
As quickly addressed in the {section \ref{AdmissibleDofs}} and as it will be shown in the numerical results, a classical construction of the space $\Hk$ implies the degeneration of some normal functions into internal ones. This is a consequence of the coordinate-wise freedom provided on the boundary from the definition of the set $A_k$. Therefore, to allow a parallel with the Raviart-Thomas elements from the lowest order on and to fulfil the optional {condition \ref{enum:4}}, one can consider replacing the boundary conditions $u|_{\partial K}\in \mathcal{H}_{l_1}(\partial K)$ in $A_k$ to obtain the reduced space
\begin{align}
	\Hk =& \{ u \in H^1(K),\,u|_{\partial K}\equiv 1,\, \Delta u \in \mathbb{Q}_{m_1}(K) \}^d \nonumber\\
	&\oplus x \,\{ u \in H^1(K),\,u|_{\partial K}\in \mathcal{H}_{l_2}(\partial K),\, \Delta u \in \mathbb{Q}_{[m_2]}(K)\}.
\end{align}
There, the coordinate-wise freedom on the boundary is reduced and the normal \Dofs can be set as in the classical Raviart-Thomas setting. Furthermore, contrarily to the general case, any definition of $l_2$, $m_1$ and $m_2$ leads to an $H(\div)$-conformity ready space.

\subsubsection*{Example of a reduced two dimensional element.}
To emphasise the parallel with the Raviart-Thomas setting on the boundary, we reduce the previous example and derive the corresponding reduced discretisation framework.
\begin{definition}[Reduced space]
	\begin{align*}
	\mathbb{H}_k(K) =& \left\{u\in H^1(K),\,u|_{\partial K} \equiv 1,\, \Delta u\in\mathbb{Q}_{k-1}(K)\right\}^2\\
	& \oplus \begin{pmatrix}x\\y\end{pmatrix}\left\{
	u\in H^1(K),\,u|_{\partial K}\in \mathcal{H}_k(\partial K),\, \Delta u\in\mathbb{Q}_{[k-1]}(K)\right\}
	\end{align*}
\end{definition}
Its dimension then naturally reads :
\begin{equation*}
\dim \mathbb{H}_k(K) = \mathfrak{n}(k+1)+2k(k-1)-\mathbbm{1}_{k>0}.
\end{equation*}
Therefore, exactly $k+1$ normal functions per edge can be designed, fitting the framework of Raviart-Thomas. As this matches the dimension of $\mathbb{Q}_k(f)$, all the freedom is required to entirely determine the global normal component. Thus, as a straightforward reduction of the general case, the $H(\div,\,K)$-conformal element presented in the {section \ref{Sec:2Dcase}} simplify to the following \Dofss.
\begin{table}[h!]
\begin{center}\resizebox{1.0\linewidth}{!}{
	\begin{tabular}{|c|c|c|c|}
		\hline
	 Core normal Moments & Misc moment & Core internal moment & Misc internal moment\\\hline
		${\Large \substack{\\ \int_{f_j}{}{x^i\,q\cdot n}{\gamma},\\ \\ \\  \\
		 	\forall j\in\irange{1}{\mathfrak{n}},\,i\in\irange{1}{k}\\ \\}}$
		&${\Large \substack{\int_{f_j}{}{q\cdot n}{\gamma}\\ \\ \\ \\ j\in\irange{1}{\mathfrak{n}}}}$
		&$\substack{\int_{K}{}{\begin{pmatrix}x^ly^m\\0			\end{pmatrix}\cdot q}{x\mathrm{d} y} \,\, \text{ and } \,\,
			\Int{K}{}{\begin{pmatrix}0\\x^my^l\end{pmatrix}\cdot q}{x\mathrm{d} y} ,\,\\
			\phantom{\int}\\
			\forall l\in\irange{0}{k},\, m\in\irange{0}{k-1}\\s.t.\,(l,\,m)\neq(k,\,k-1)}$ &${ \substack{\int_{K}{}{\begin{pmatrix}x^{k}y^{k-1}\\x^{k-1}y^{k}		 \end{pmatrix}\cdot q}{x\mathrm{d} y}\\ \\ \\ \\ \\ \\}}$\\\hline
	\end{tabular}}
\end{center}\caption{Degrees of freedom of the element \IIa \space defined within the reduced setting}
\end{table}
 Note that here too,  for defining the projections \emph{\eqref{eq:DofsIntHk2D2}} one can work with any basis of $\mathcal{P}$ instead of the presented canonical basis.

	\section{Numerical results.}
\label{NumRes}
We explore the properties of the main element \IIa{} and its variant \IIb{} presented in the previous section by investigating their basis functions, for both the general framework and the reduced one.
We particularly focus on the normal component $\varphi\cdot n|_{\partial K}$ of representative basis functions $\varphi$ on the boundary of the element $K$. Those basis functions have been constructed by tuning a natural basis of the space $\mathbb{H}_k(K)$ towards the selected sets of \Dofs through a transfer matrix.\newline

As an example, we consider the non-convex nine-edges polygon presented in the {figure \ref{poly}} on which the elements are built. In all the results, the polynomial projectors used in the definition of the \Dofs were chosen as Hermite polynomials:  experimentally we have observed that this improves the conditioning of the linear system.

\subsection{General setting.}
We start by considering the spaces and elements described in the {section \ref{Sec:2Dcase1}}.
\begin{figure}[h]
	\centering
	\includegraphics[width = 0.4\textwidth]{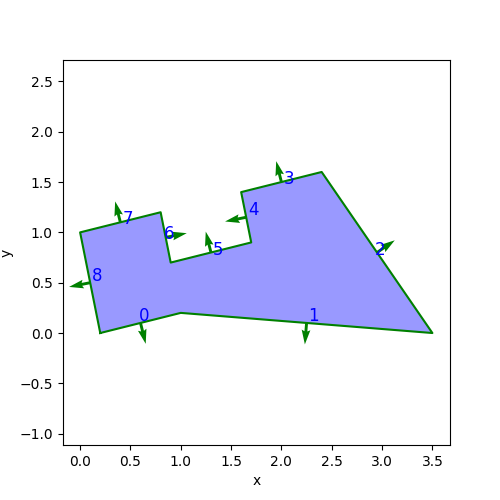}
	\includegraphics[width = 0.55\textwidth]{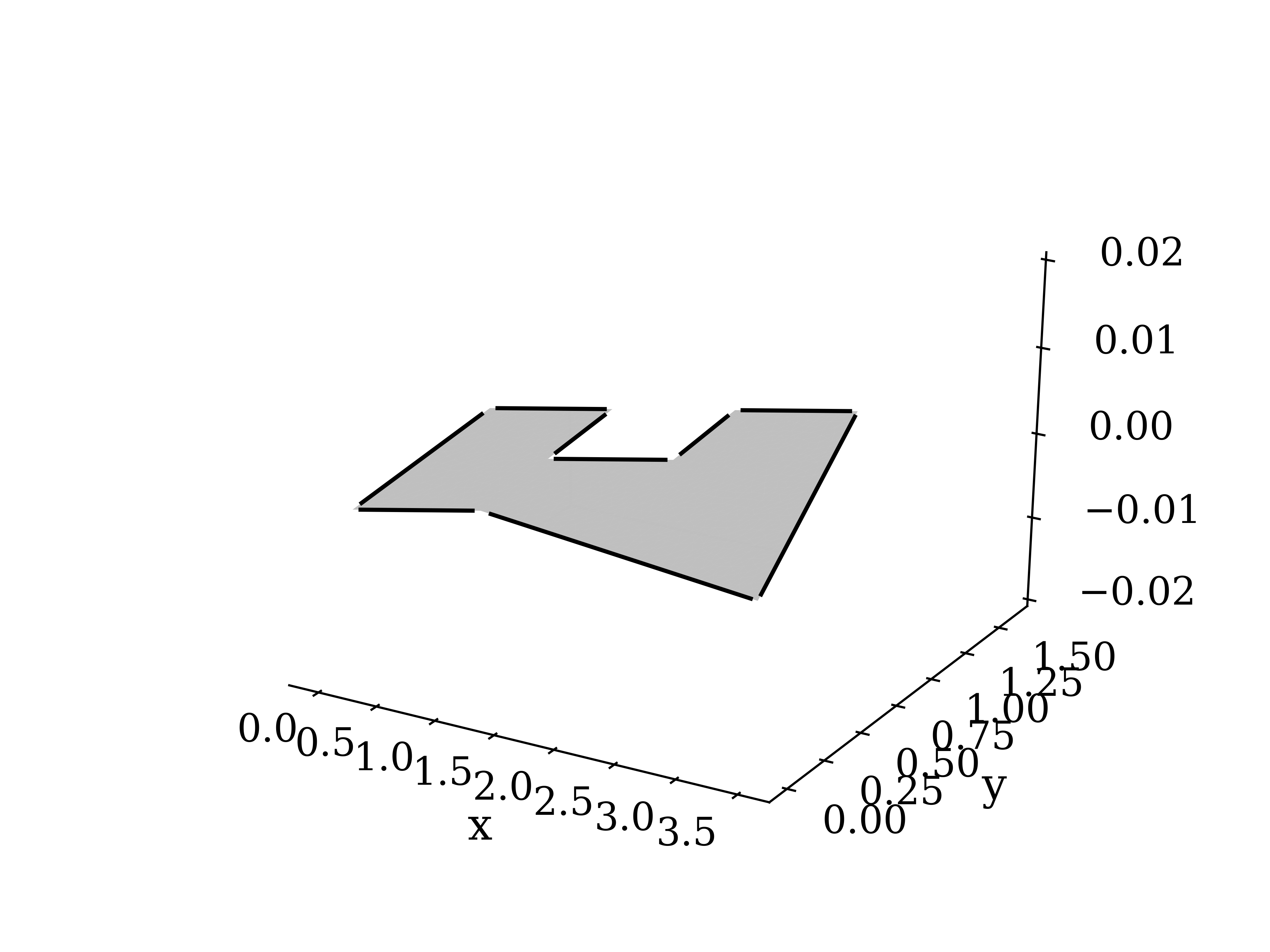}

	\caption{\label{poly}Left: Considered polygon. Right: Normal component of a representative internal basis functions plotted on every edge.}
\end{figure}
First of all, we have investigated the behaviour of the internal basis functions. The normal component of them is shown in the right of figure \ref{poly}. As wished, the basis functions corresponding to internal degrees of freedom vanish on the boundary. This can been seen on the right figure where the function is plotted in the plane $z = 0$.

\begin{figure}[h!]\centering
	\includegraphics[width = 0.4\textwidth]{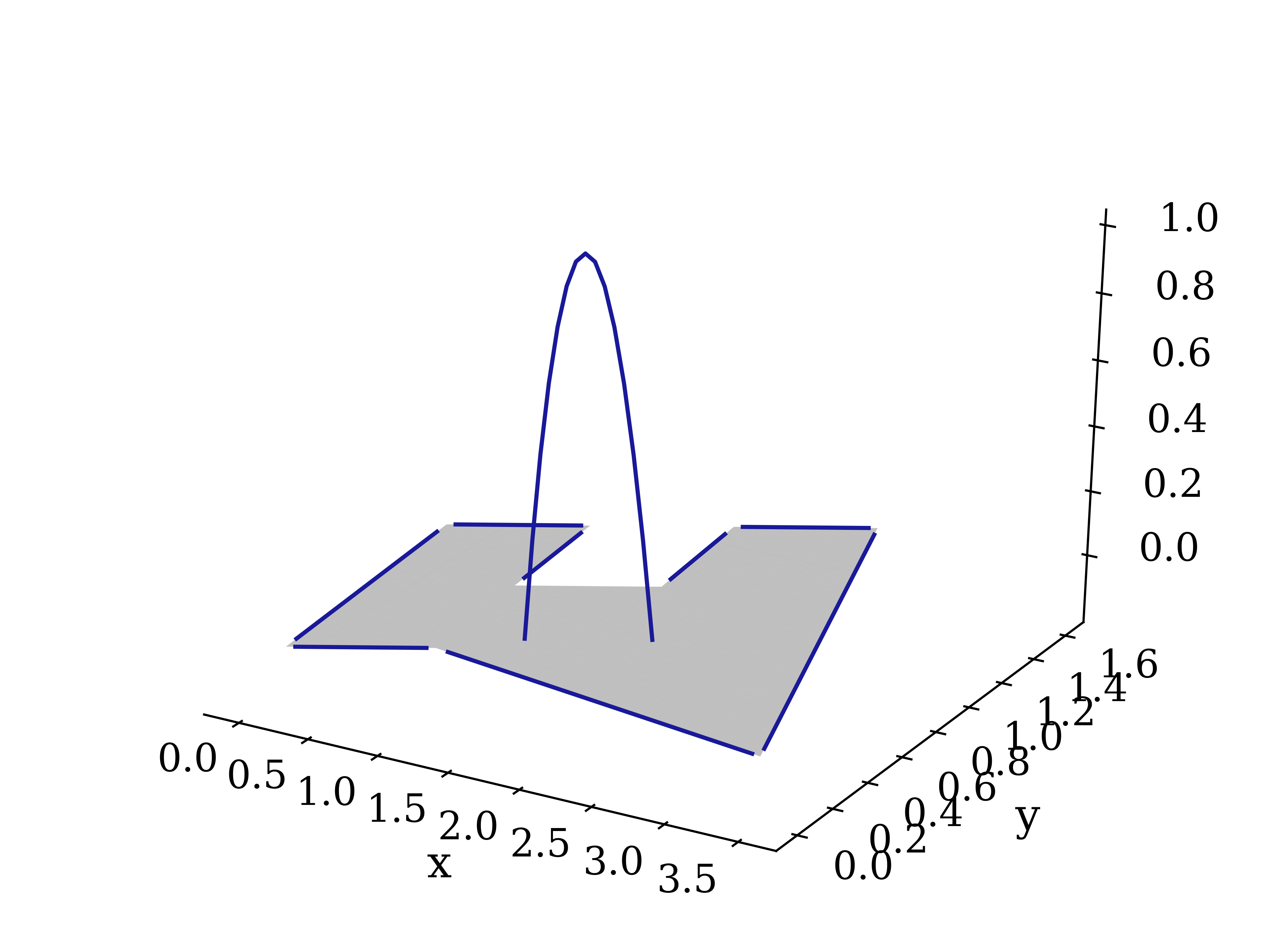}
	\includegraphics[width = 0.28\textwidth]{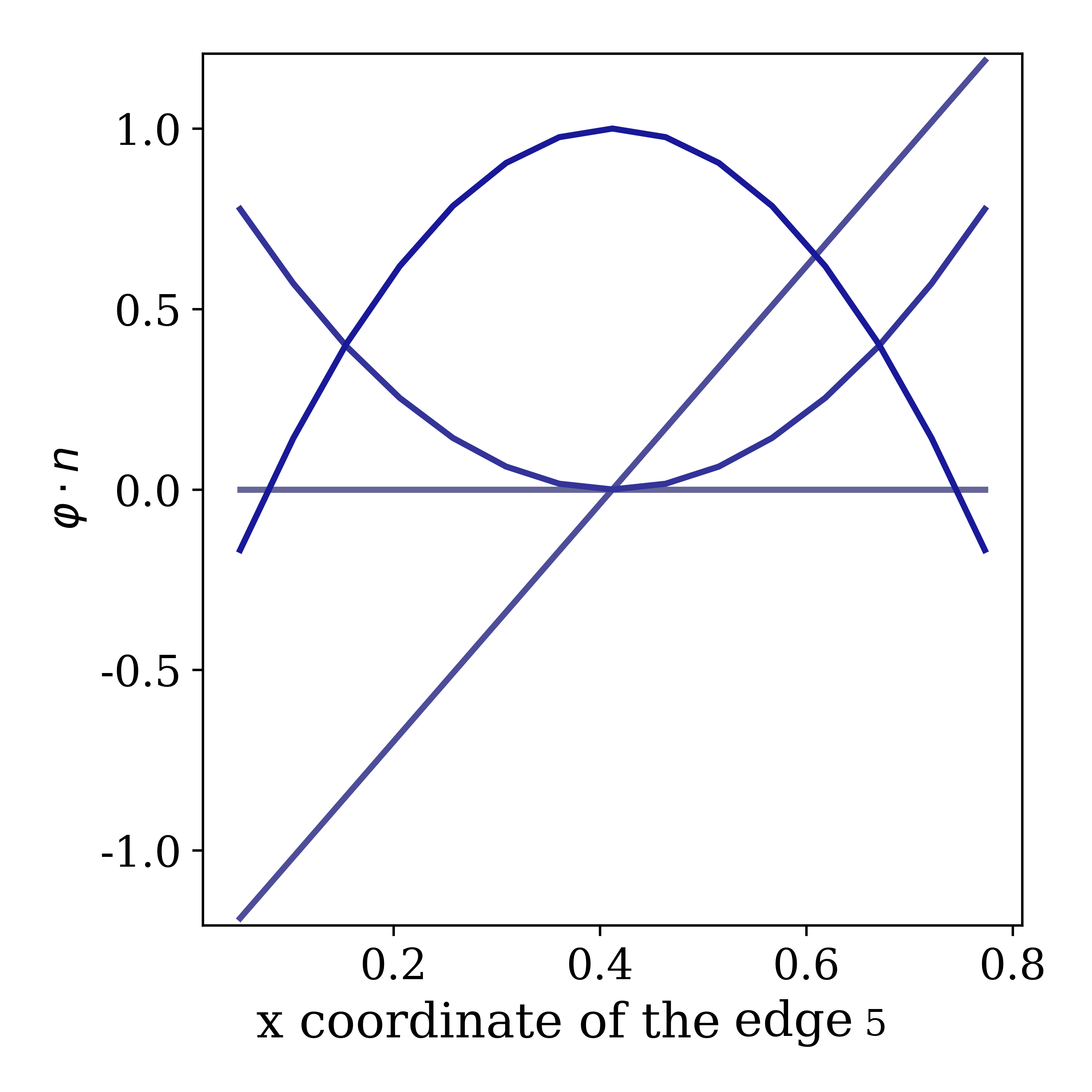}
	\includegraphics[width = 0.28\textwidth]{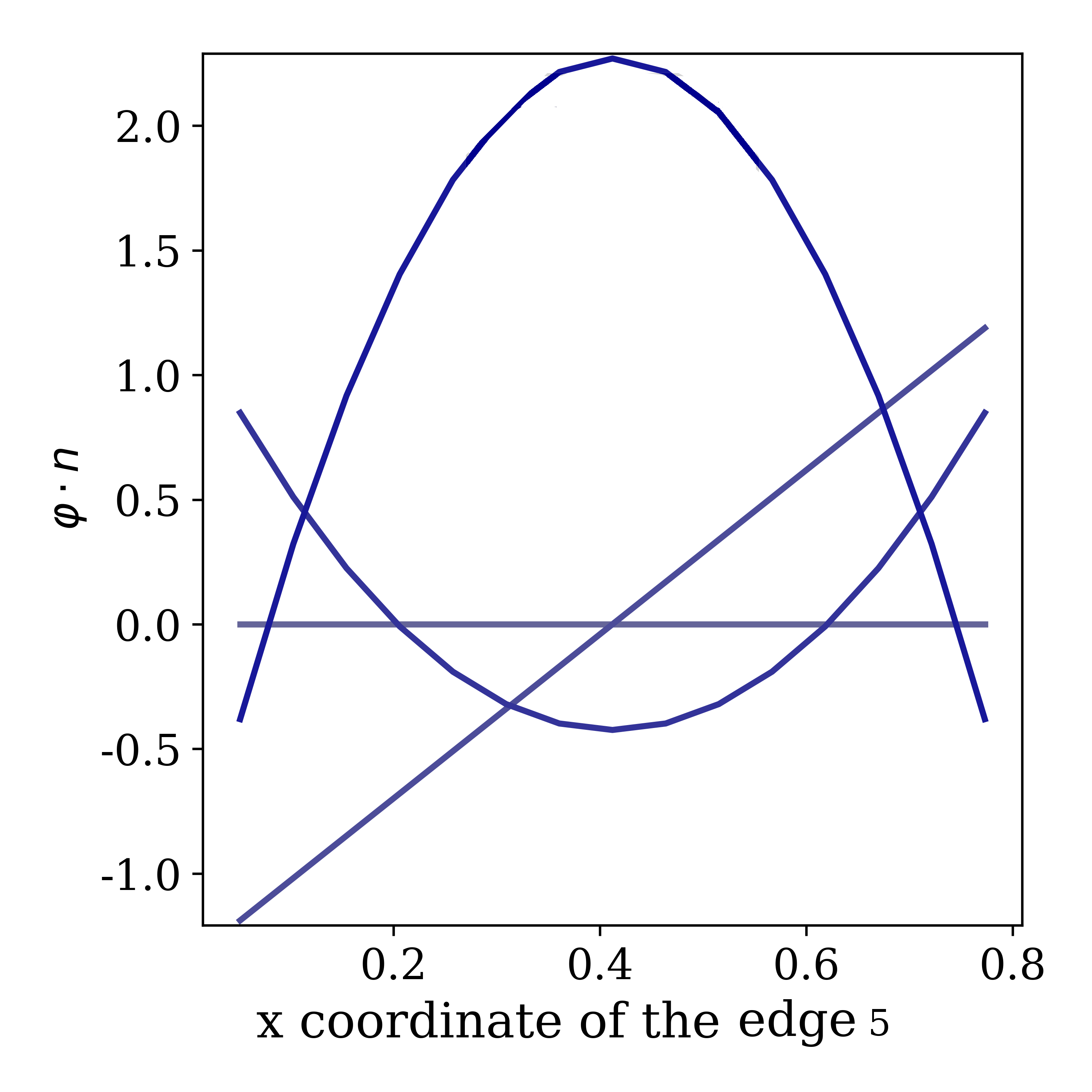}

	\caption{Left: normal component of one representative of the normal basis functions for the element \IIb{} and $k=2$
	along the edges. Middle: normal component of all the functions generated from the edge number $5$, plotted on the edge number $5$. Right: as middle, for the element \IIa.\label{IIbnorm}}
\end{figure}

In order to study the behaviour of the normal basis functions on the boundaries, we have considered the case $k=2$ where we expect five basis function to have a quadratic normal component. We have plotted in the left most  side of the figure \ref{IIbnorm} the normal component of one of the normal basis functions associated to the element \IIb. One can observe that its support is contained on one single edge.

\begin{figure}[h!]\centering
	\includegraphics[width = 0.45\textwidth]{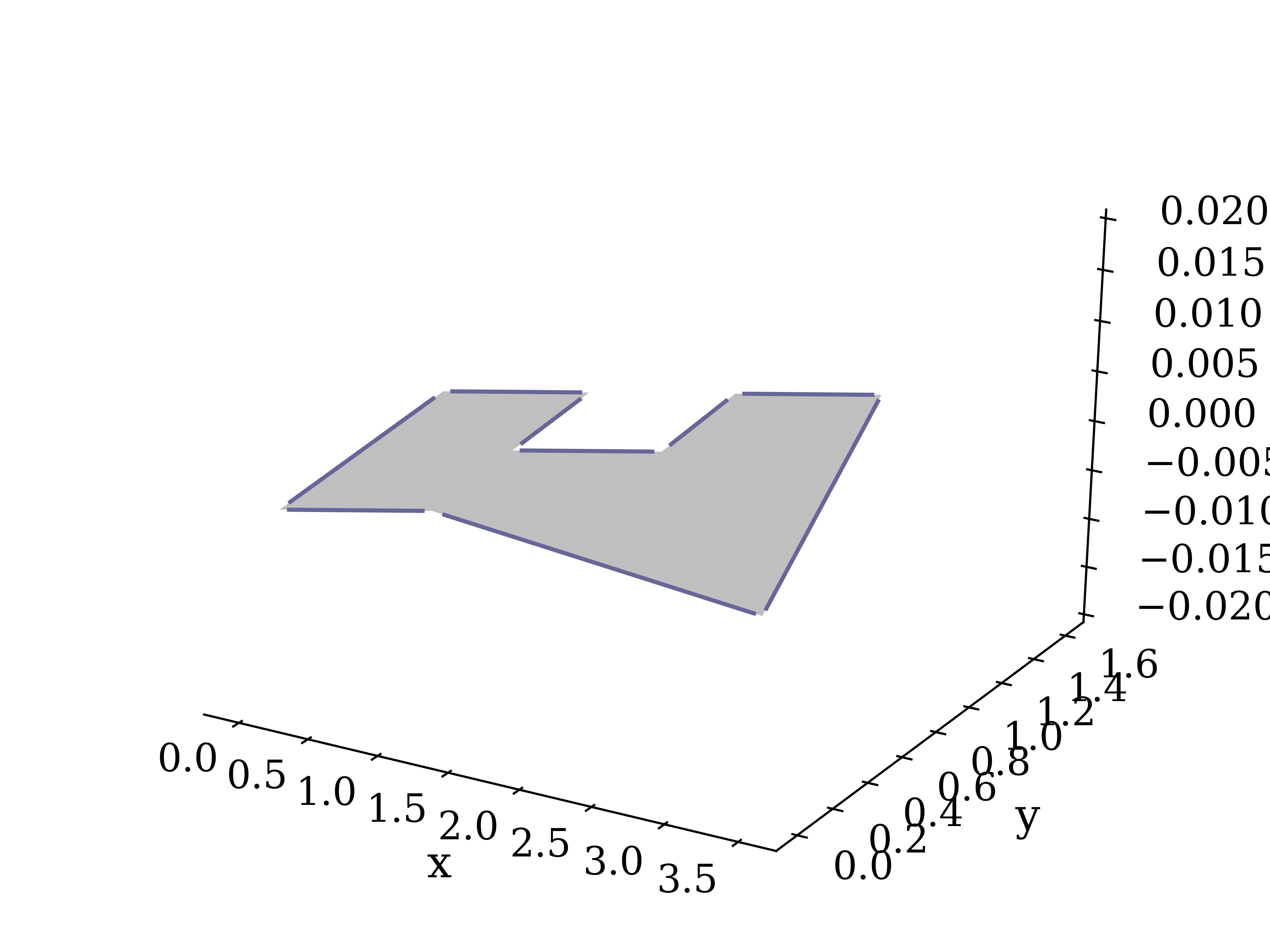}
	\includegraphics[width = 0.45\textwidth]{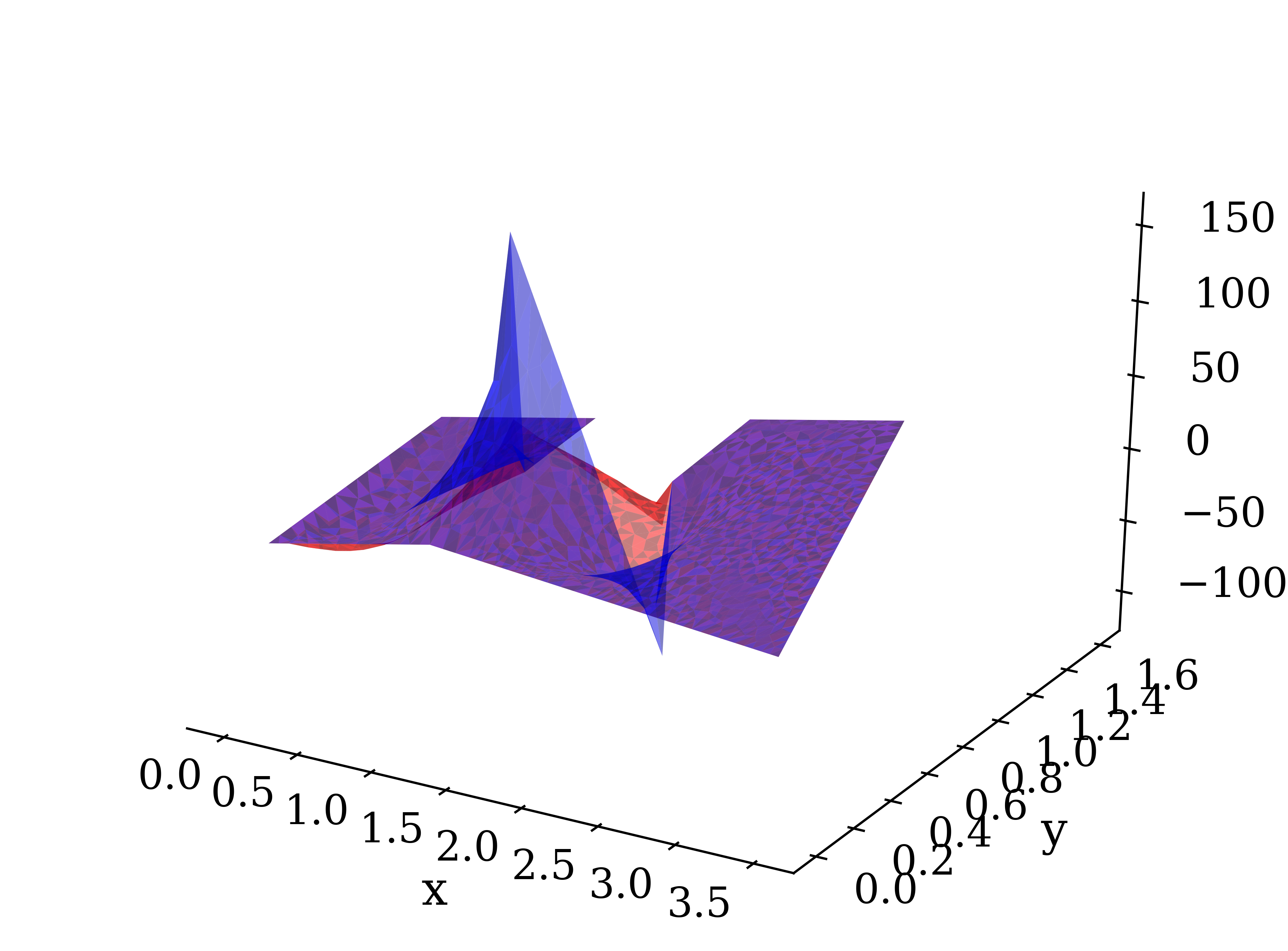}
	\caption{Degeneration of a degenerating normal basis functions' representative in the case $k=2$, for the element \IIb. Left: normal component on all the boundaries. Right: internal behaviour of the basis functions. \label{IIbvanish}}
\end{figure}

We then have plotted all the basis functions associated to the edge number $5$ on the figure \ref{IIbnorm}, for both the configurations \IIb{} (middle) and \IIa{} (right). One can first notice that their normal components, plotted in the middle graph of the figure \ref{IIbnorm},  are polynomial of degree $k\leq2$, that together generate the space $\mathbb{P}_2(\mathbb{R})$. Observing further, it appears that the normal component of two basis are vanishing, that is $\varphi\cdot n =\varphi_1\cdot n_{5,\,x} +\varphi_2\cdot n_{5,\,y}=0$. This generates this straight line equal to zero in the graph. Indeed, those two basis functions characterise the coordinate-wise freedom $\varphi_1\cdot n_{5,\,x}\neq 0$ and $\varphi_2\cdot n_{5,\,y}\neq 0$. This additional freedom is  not reflected through the global term $\varphi\cdot n$ as addressed in the example \ref{ex:Struct} and in the remark \ref{reducing}.  
This comes from the fact that only three basis functions are required to generate $\mathbb{P}_2(\mathbb{R})$, where the global component $\varphi\cdot n$ lives. The two components $\varphi_1n_{5,\,x}$ and $\varphi_2n_{5,\,y}$ of the vector $\varphi \circ n$ are compensating themselves. Those functions are nevertheless regular within the polygon $K$ and not identically vanishing on $K$ (see figure \ref{IIbvanish}, where the left graph represents the value of the normal component $\varphi\cdot n$ on the boundary and where the right graph represents the components $\varphi_1$ and $\varphi_2$  on the element $K$). They can therefore be reclassified as internal basis functions. Note that this can be suppressed when using the reduced setting, as one can observe below in the section \ref{reducednum}.

\begin{figure}[h!]\centering
	\includegraphics[width = 0.28\textwidth]{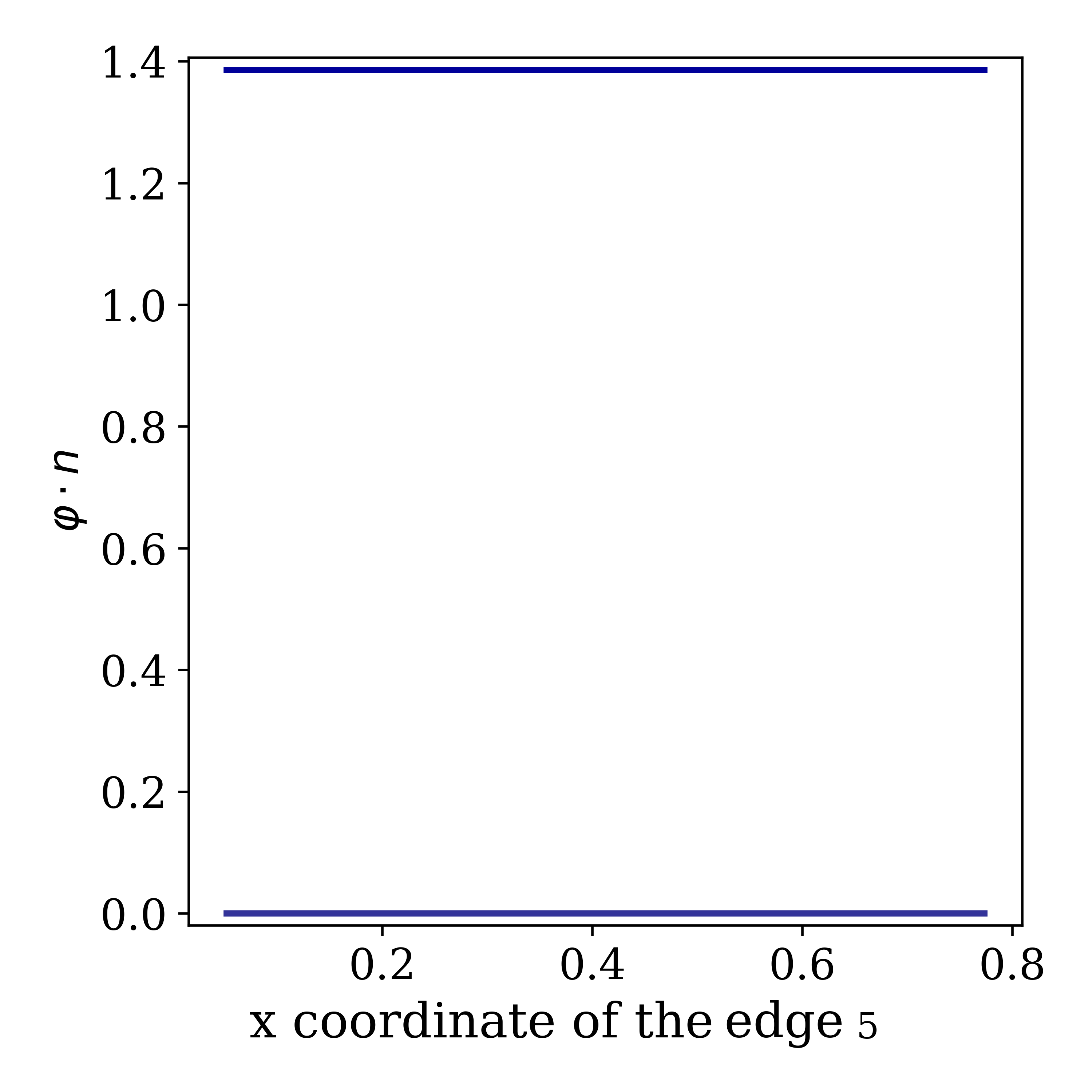}\quad\quad
	\includegraphics[width = 0.28\textwidth]{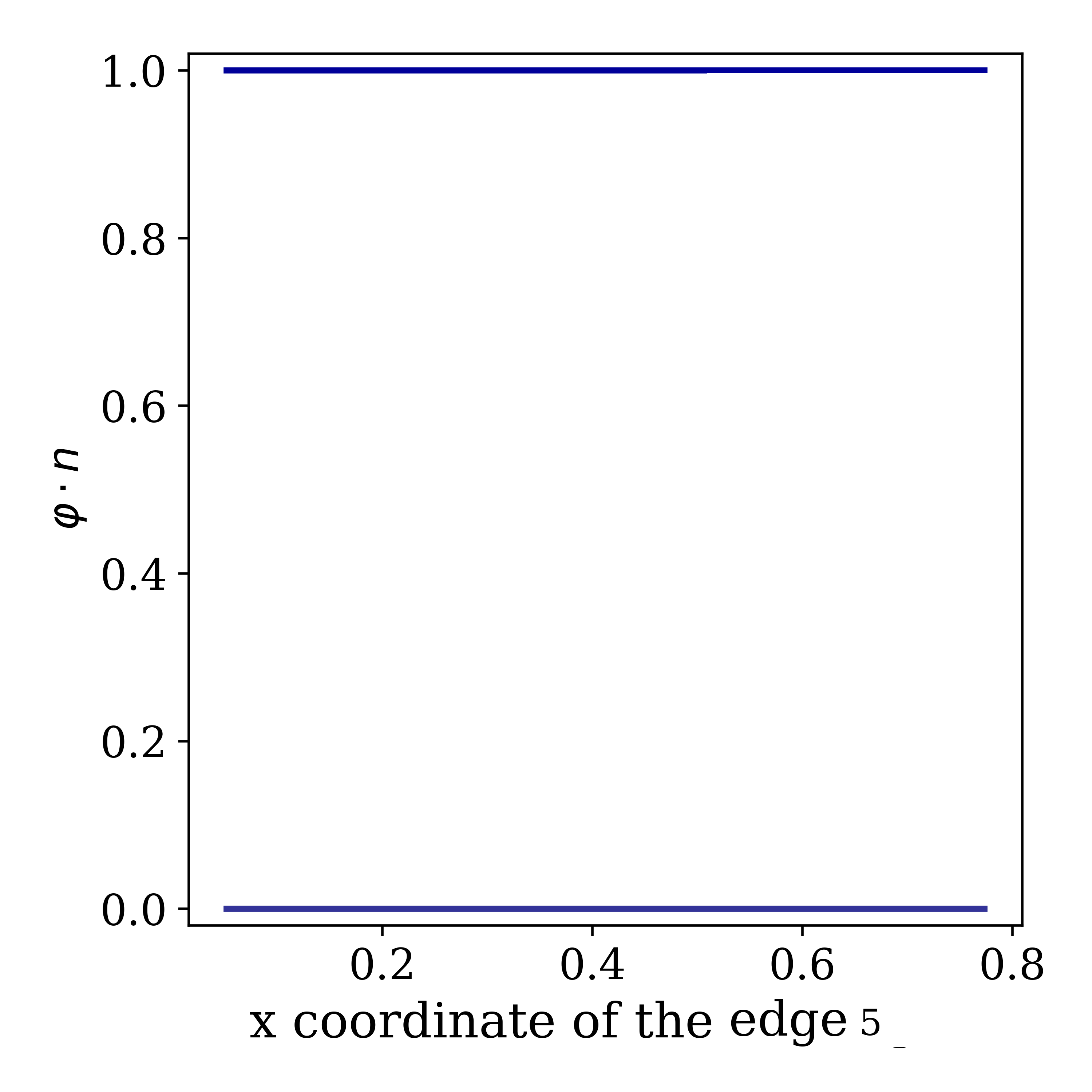}
	\caption{Scaling of the non-vanishing basis functions generated from the edge number $5$ when $k=0$. From left to right: \IIa{} and \IIb. \label{scaling}}
\end{figure}

Finaly, one can consider the scaling of the basis functions by plotting the normal basis functions corresponding to the configurations \IIa{} and \IIb{} in the lowest order case, \emph{i.e.} for $k=0$ (see figure \ref{scaling}). There, only the configuration \IIb{}-using a point-wise value- scales to one. The fully moment-based configuration \IIa{} scales to another constant that depends on the edge's length and orientation with respect to the axes. This example emphasises that the configuration \IIb{} leads to basis functions which share similar properties like the Raviart-Thomas elements.

\subsection{Reduced setting.}
\label{reducednum}
As a last example, we derive some results obtained for the reduced element \IIb, offering a complete parallel with the Raviart-Thomas elements on the boundary by suppressing the further coordinate-wise liberty provided by the general setting. The internal basis functions being unchanged from the general setting, they are not represented.

\begin{figure}[h!]\centering
	\includegraphics[width = 0.45\textwidth]{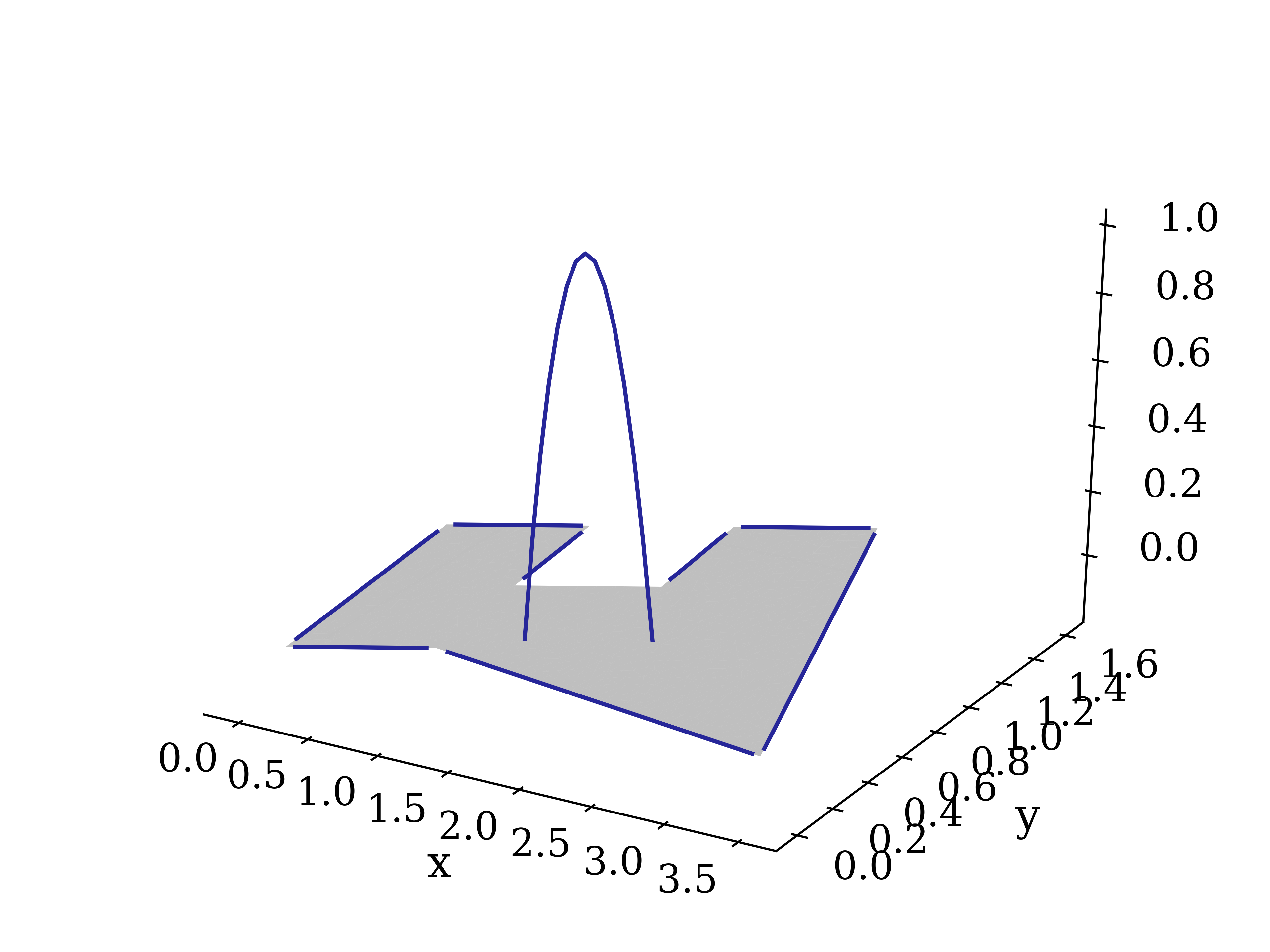}
	\includegraphics[width = 0.45\textwidth]{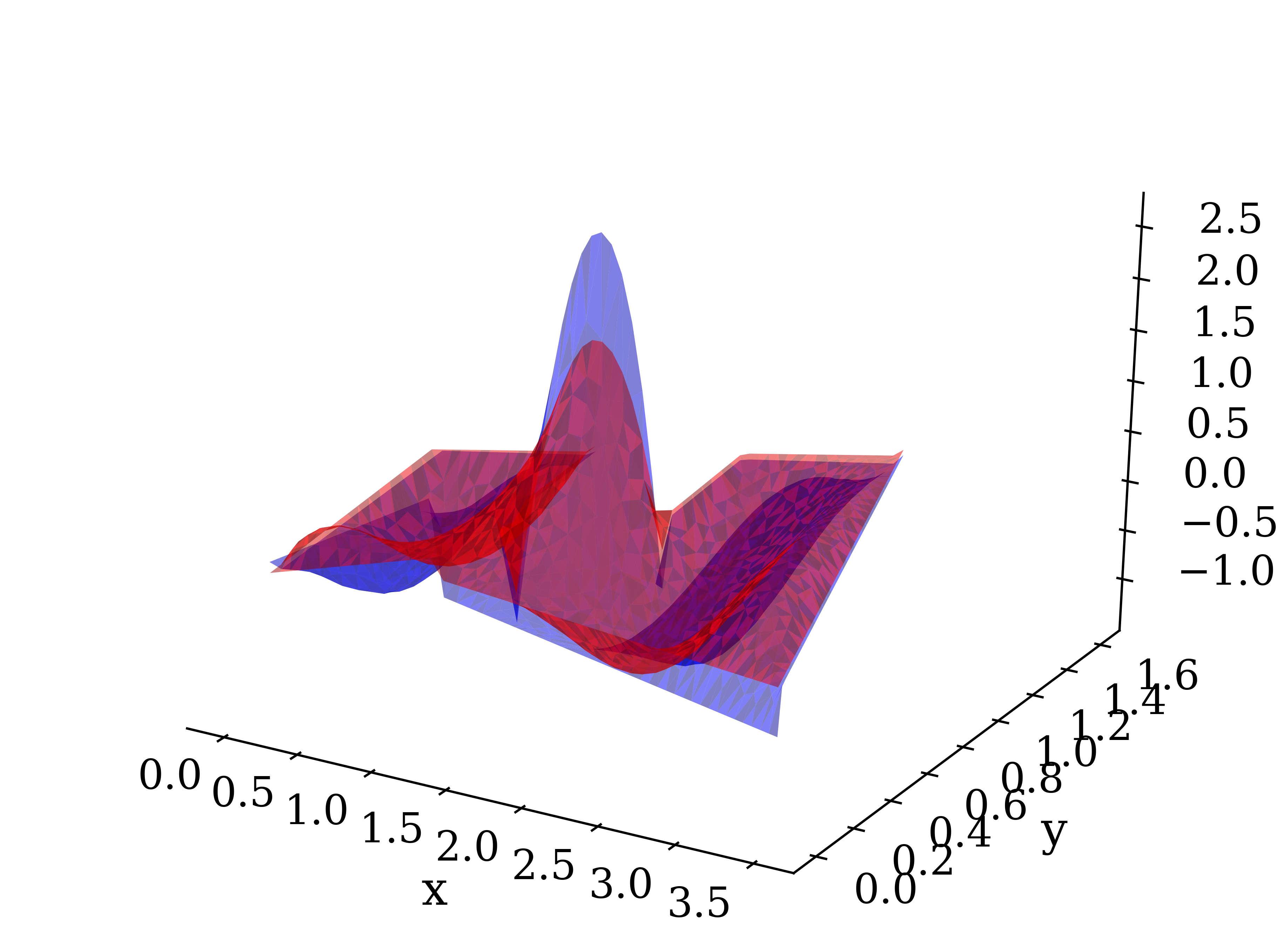}

	\includegraphics[width = 0.28\textwidth]{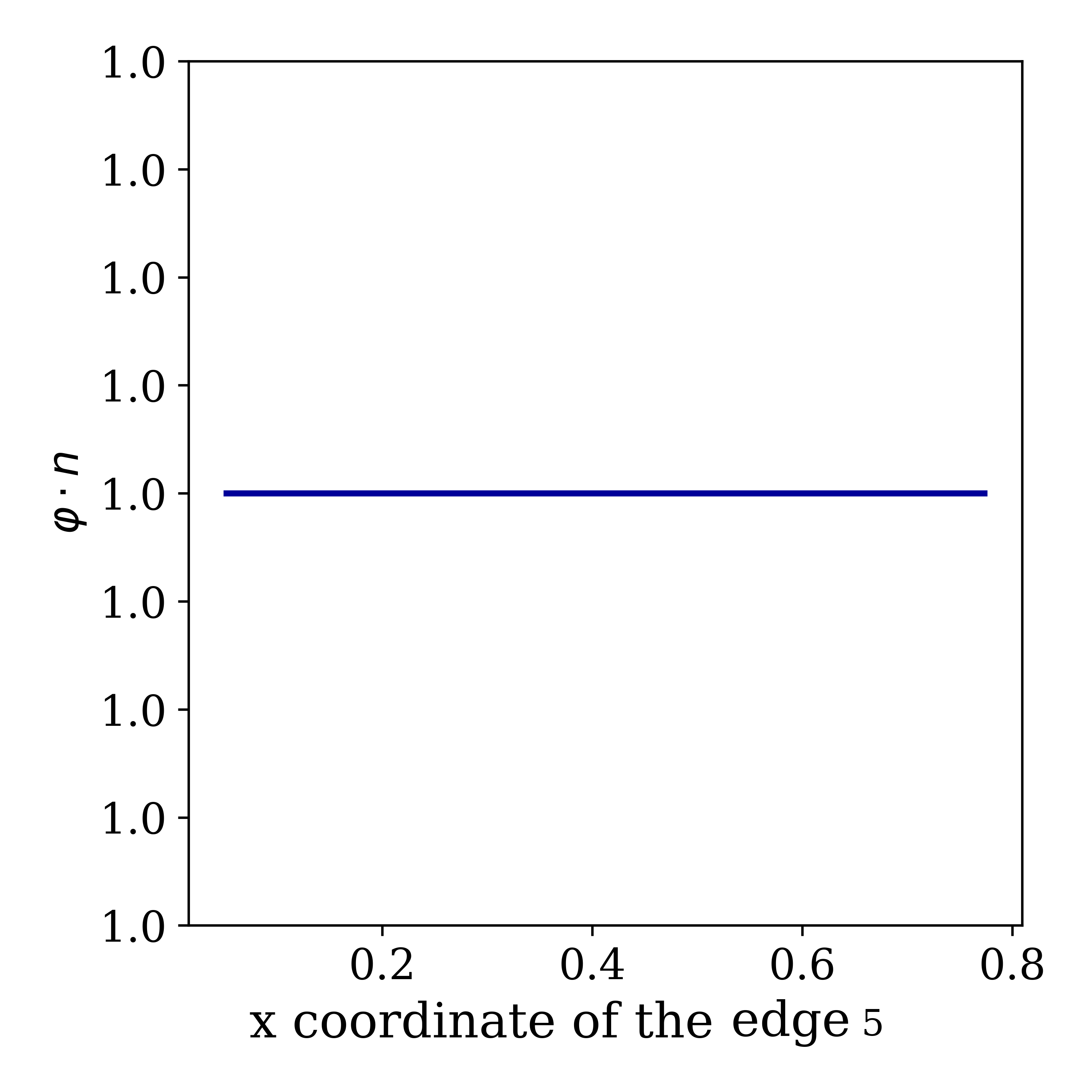}
	\includegraphics[width = 0.28\textwidth]{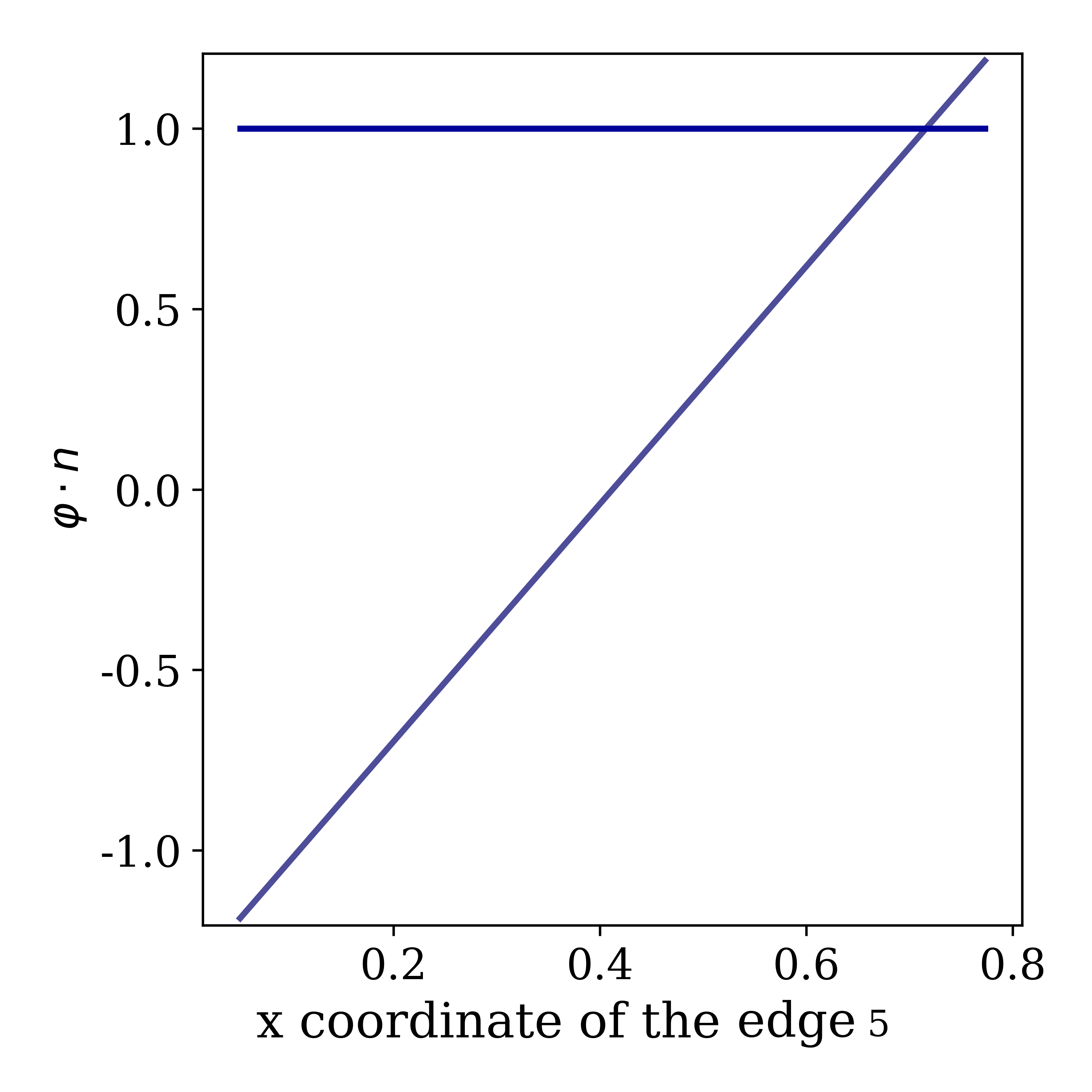}
	\includegraphics[width = 0.28\textwidth]{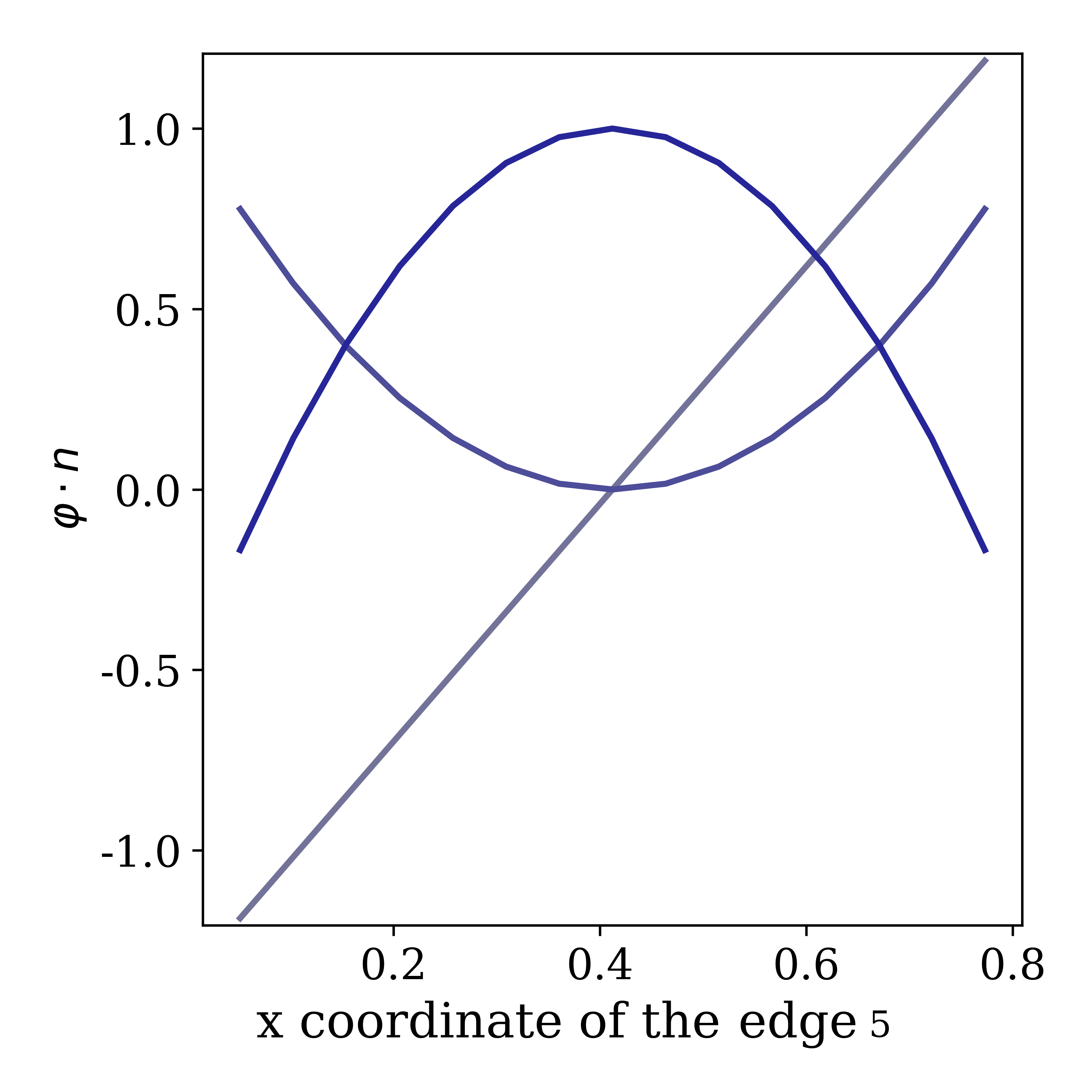}
	\caption{\label{reducedresults}Top left: regularity of the components of one representative of the basis functions for the reduced element \IIb{} and $k=2$ within the element. Top right: its normal component along the boundary. Bottom: normal component of all the functions generated from the edge number $5$, plot on the edge number $5$. From left to right: $k=0$, $k=1$, $k=2$.}
\end{figure}

Indeed, one can observe on the bottom of the {figure \ref{reducedresults}} that there is no more degenerating normal basis functions. Therefore, all normal basis functions are acting globally to characterise the polynomial behaviour of functions of the reduced $\Hk$ space on the boundary. Furthermore, one can observe that the scaling of the basis functions corresponding to lowest order element, as well as the amplitude of the basis functions describing the higher order ones, make the discretisation framework reliable.
\begin{remark}
	The shape of the normal component of the basis functions is driven by the definition of the projectors $p$ in the normal \Dofss. Changing the basis of the projectors then allows to enforce wished shape of the basis functions of $\Hk$ while keeping the regularity and order of the discretisation. Shifting them by modulating the offset directly from the definition of the \Dofs to enforce their positivity is equally possible.
\end{remark}


	\section{Conclusion}

Motivated by defining a flux reconstruction scheme on general polytopes \cite{abgrall:hal-01820176},
we have  developed a new $H(\div)$-conformal discretisation framework that can be set up on any polytope,
not necessarily convex. It merges the flexibility of the Virtual Element setting with the properties of the Raviart-Thomas elements on the boundaries.

The introduced finite dimensional spaces are vectorial and allow a lot of flexibility in the definition of the \Dofss. In particular, the choices of discretisation quality and \Dofs on the boundary are independent from the ones made within the element.


The discretised quantities benefit from an extensive coordinate-wise freedom. Therefore, upon the choice made while selecting the
\Dofss, some dual normal basis functions may be reclassified into internal ones. Thus, to allow a complete parallel with the Raviart-Thomas setting on
the boundary from the lowest order on, one may construct straightforwardly a reduced space, along with reduced elements.

Last, we detailed a particular example of a discretisation framework through a series of spaces and the definition of a particular element.
It could be observed that in both general and reduced frameworks, the type of \Dofs (point-wise values or moments) impacts the scaling of the dual basis functions.
This can typically be observed in the lowest order case of the given examples, where only the dual basis functions of the element \IIb{} scale to one.

An important topic for further research is the exploration of projectors from the introduced spaces onto polynomial ones, in a way similar to those already constructed in \cite{da2016h} and used in the Virtual Elements Methods. Especially, a suitable extension of those projectors to the spaces introduced here may be inferred from a close investigation of those projectors. Once the projectors have been defined, we can apply those discretisation spaces in a more practical context, as by example employing the introduced spaces in a finite element framework.

To conclude, let us point out again that the results presented here are already useful from a theoretical point of view. Indeed, they first guarantee that the considerations about FR schemes on general polytopes hold, and guarantee that the conjecture about the correction functions made in \cite{abgrall:hal-01820176} is correct. Secondly, it opens the door to a more general framework in context of FE, direction that will be further investigated in the future.


	\section*{Acknowledgements.}
	E. Le M\'el\'edo and P. \"Offner have been funded by SNF project 200020\_175784 "Solving advection dominated problems with high order schemes with polygonal meshes: application to compressible and incompressible flow problems".
	P. \"Offner has also been  funded by an University of Z\"urich Forschungskredit grant. The constructive criticisms of the three referees are also acknowledged.
	
	\appendix
	\section{ A note on further possible configurations}
\label{furtherconfig}

As example, we only detail in this paper two declinations of one possible configuration of \Dofss. There exists many more possibilities,  and their choice impact the properties of the elements. In particular, it is possible to focus on a coordinate-wise boundary characterisation of the quantities living in a general space $\Hk$, rather than the global focus presented here. In the general setting, this choice yields a degeneracy of only one normal basis function, thus moving a bit away from the Raviart-Thomas spirit. For interested readers, a presentation of this possibility for both the general and the reduced space is available in \cite{abgrall2019class}, along with further investigations on various configurations.

Note also that different choices also have a very strong impact on the conditionning number of the linear systems to solve. The solution we have shown in this paper is the one that offers the best compromise. It is also the one that is the closest from the classical RT framework.

	\section{Proofs}
\label{propro}
\begin{proof}[Proof {of  Proposition \ref{Prop:Hdiv}}]
We start by deriving the first statement. By construction, any $q\in\Hk$ can be decomposed into $q = q_0 + x\,q_1 $ for some $q_0\in (A_k)^d$ and $q_1\in B_k$.
	Therefore, on the boundary of $K$ one has $q\,\cdot n |_{\partial K} = q_0\cdot n|_{\partial K} + (x q_1\cdot n)|_{\partial K}$. As the functions $q_1$ is scalar, this quantity can also read $q\,\cdot n |_{\partial K} = q_0\cdot n|_{\partial K} + q_1(x\cdot n )|_{\partial K}$ {by linearity and commutativity of the dot product.}

	Since for every face $f$ of $K$ the term $x\cdot n|_{f}$ is constant, it reduces to $q\,\cdot n|_{f} = q_0\cdot n|_{f} + c_f\,q_1|_{f}$ on each face $f$ for a constant $c_f\in\R$ depending only on the face layout and position with respect to the axes and origin. Therefore, since $q_0|_{f}\in(\mathbb{Q}_{l_1}(f))^d$ and $q_1|_{f}\in\mathbb{Q}_{l_2}(f)$, $q\,\cdot n |_{f} \in \mathbb{Q}_{\max \{l_1,\,l_2\}}(f)$. And since it is valid for any face $f\in\partial K$, we finally get that $q\,\cdot n|_{\partial K} \in \mathbb{H}_{\max\{l_1,\,l_2\} }(\partial K)$

 Let us now derive the divergence property within the cell. Any $u\in\Hk$ can be written under the form
$u = \tilde{q} + x\,q$
for some functions $q\in H^1(K)$ and $\tilde{q}=(\tilde{q_1},\,\cdots,\,\tilde{q_d})^T\in (H^1(K))^d$ such that
\begin{equation}
\label{eq:Propq}
\left\lbrace \begin{aligned}
\Delta q \in \mathbb{Q}_{[m_2]}(K)\\
q|_{\partial K}\in\mathbb{Q}_{l_2}(\partial K)
\end{aligned}\right.
\quad \text{ and } \quad
\left\lbrace \begin{aligned}
\Delta \tilde{q_i} \in \mathbb{Q}_{m_1}(K)\\
\tilde{q_i}|_{\partial K}\in\mathbb{Q}_{l_1}(\partial K),
\end{aligned}
\right.\quad {\forall i\in\irange{1}{d}.}
\end{equation}
We have
\begin{align*}
\div(u) &= \sum\limits_{i=1}^{d}\partial_{x_i}(x_i\,q)+\sum\limits_{i=1}^{d}\partial_{x_i}\tilde{q_i}= \sum\limits_{i=1}^{d}\left(q+x_i\partial_{x_i}q\right)+\sum\limits_{i=1}^{d}\partial_{x_i}\tilde{q_i}\\
&=\underbrace{d\,q}_{\in L^2(K)} + \sum\limits_{i=1}^{d}\big(x_i\underbrace{\partial_{x_i} q}_{\in L^2(K)}\big) + \sum\limits_{i=1}^{d}\underbrace{\partial_{x_i} \tilde{q_i}.}_{\in L^2(K)}
\end{align*}
Since by \eqref{eq:Propq} we have $\nabla\cdot q\in L^2(K)$, it comes that for any $i\in\irange{1}{d}$;
$x_i\,\partial_{x_i} q\in L^2_{\mathrm{loc}}(K)$.
As $K$ is compact and bounded, we have $L^2_{\mathrm{loc}}(K) = L^2(K)$ and $\div \,q\in L^2(K).$ As a by-product, note that we can derive
$
\nabla\cdot(x\nabla q)	 = \nabla\cdot q+x\Delta q$, 
where $\Delta q\in \mathbb{Q}_{\max{[m_1,\,m_2+1]}}$ and $x\Delta q\in \mathcal{C}^{\infty}(K)$.
\end{proof}

\begin{lemma}
	\label{LemmaImportant}
	The configurations \IIa{} and \IIb{} are sets of \Dofs leading to unisolvent elements when endowed in $\Hk$.
\end{lemma}

\begin{proof}[Proof { of the lemma \ref{LemmaImportant}}]
	\label{proofimportant}
	\emph{We refer to the functions $p_k$ by the term \enquote{kernel}, while using the term \enquote{integrand} to represent the term $q\cdot p_k$. Immediate transfer of this designation apply to the normal moment based \Dofss.}

	We first sketch the proof. To begin with, let us point out that the key lies in the {assumptions \ref{Ass:Dofs0}}  ensuring the linear independence of the set of point-wise values and moment's integrands. The linearity of the integral operators transfers then this independence to the moments themselves, characterising any function of $\Hk$ independently on the boundary and within the cell. We proceed in three steps.
	\begin{enumerate}
		\item First, we show that the internal characterisation of the function does not impact the normal one, allowing the determination to be done distinctively within the element and on the boundary.

		\item Then,  we show that selecting the appropriate number of \Dofs in any of the sets \IIa{} or \IIb{} ensures a unique characterisation on the boundary. We use the fact that the kernels are scalar polynomials while the functions of $\Hk$ are vector polynomials.

		\item Lastly, we consider the interior of the element where the characterisation is done through projections over linearly independent sets. Those projections of functions in $\Hk$ are indeed neither identically null nor identically identical (\emph{i.e.} they differ at least on a subset of non-zero measure).
	\end{enumerate}
	Let us detail this determination process more in details.

\noindent \textit{Step 1.} Let us first recall that the space $\Hk$ is built from blocks of independent functions. In particular, the boundary behaviour of functions living in $\Hk$ is independent of their behaviour within the inner cell. Therefore, by the structure of $\Hk$ and making use of the superposition theorem, any function $q\in\Hk$ reads
	\begin{equation}
	\label{eq:qspli}
	q = f\mathbbm{1}_{\partial K} + g\mathbbm{1}_{\mathring K}
	\end{equation}
	for two functions $f$ and $g$ belonging to $\Hk$. As a consequence, characterising a function $q\in\Hk$ comes down to characterising the independent functions $f$ and $g$ on the distinct supports $\partial K$ and $\mathring K$, respectively.
	Note also that necessarily, $f|_{f_j}\in\bigtimes_{i=1}^d\mathbb{Q}_{\max\{l_1,\,l_2+1\}}(\R^{d-1})$ for any face $f_j\in\partial K$. We show that any admissible extraction (in the sense of the {admissibility conditions \ref{Ass:Dofs0}}) from either of the two sets of \Dofs (\IIa, \text{internal}), (\IIb,\, \text{internal}) fully characterises the functions $f$ and $g$, independently. In all the following, the notation \IIa{} or \IIb{} refers to the corresponding set of normal \Dofs while "internal" refers to the set $\eqref{eq:DofsIntHk}$ and is identical to any of the two configurations under consideration.

	We first show that any above defined set of \Dofs preserve the independence of the boundary and inner characterisations. To this aim, we combine the relation \eqref{eq:qspli} with the all possible definitions of the \Dofss. It comes that all global normal moments lead to an expression of the form
	\[\sigma(q) = \SInt{f_j}{}{q\cdot n\, p_k} = \SInt{f_j}{}{(f\mathbbm{1}_{f_j} + g\mathbbm{1}_{\mathring{K}})\cdot n \, p_k} = \SInt{f_j}{}{f\cdot n \, p_k} \]
	for some polynomial function $p_k$ living on $\partial K$. On the other side, as $x_{jm}\in f_j$, the global \Dofs that are built from point-wise values read
	\[\sigma(q) = q(x_{jm})\cdot n  = f(x_{jm})\cdot n\,\mathbbm{1}_{f_j}(x_{jm}) + g(x_{jm})\cdot n\,\mathbbm{1}_{\mathring{K}}(x_{jm}) = f(x_{jm})\cdot n.\]
	Similar relations for coordinate - wise \Dofs can be derived, that is;
	\begin{align*}\sigma(q) &= \SInt{f_j}{}{q_{x_i} n_{x_i}\, p_k} = \SInt{f_j}{}{(f_{x_i}\mathbbm{1}_{f_j} + g_{x_i}\mathbbm{1}_{\mathring{K}}) n_{x_i} \, p_k} = \SInt{f_j}{}{f_{x_i}n_{x_i} \, p_k} \\
	\text{and \quad} \sigma(q) &= q_{x_i}(x_{jm}) n_{x_i}  = f_{x_i}(x_{jm}) n_{x_i}\mathbbm{1}_{f_j}(x_{jm}) + g_{x_i}(x_{jm}) n_{x_i}\mathbbm{1}_{\mathring{K}}(x_{jm})\\
	& = f_{x_i}(x_{jm}) n_{x_i},
	\end{align*}
	where here the terms $f_{x_i}$ simply represent the $i-\mathrm{th}$ component of the function $f$.
	Therefore, in any of the configurations \IIa{} and \IIb{} no contribution of the function $g$ representing the inner part of the cell is involved in the normal \Dofss. The mirror case is obtained with the internal moments, leading via \eqref{eq:qspli} to
	\[\sigma(q) = \SInt{K}{}{q\cdot \,p_k} = \SInt{K}{}{(f\mathbbm{1}_{\partial K} + g\mathbbm{1}_{\mathring{K}})\cdot \, p_k} = \SInt{ K}{}{g\cdot \, p_k}, \]
	where $p_k$ stands for any Poisson's function living in $\Hk$ or any polynomial function defining the second member of a Poisson's problem involved in the definition of $\Hk$. There, the function $f$ representing the boundary part of the function $q$ is not involved, that for any definition of the space $\mathcal{P}_k$ generating the internal moments. Thus, by linearity we can decompose the \Dofs $\{q\mapsto\sigma_{i}(q)\}_i$ in the following matrix.
	\begin{center}
		\resizebox{0.95\textwidth}{!}{
			\begin{tikzpicture}[decoration=brace]
			\matrix (b) [matrix of math nodes,left delimiter=(,right delimiter={)}] {
				\sigma_1\phantom{+1}\\
				\vdots\\
				\sigma_{N_N\phantom{+1}}\\
				\sigma_{N_N+1}\\
				\vdots\\
				\sigma_{N_I\phantom{+1}}\\
			};
			\begin{scope}[shift={(5.6,0)}]
			\matrix (m) [matrix of math nodes,ampersand replacement=\&,left delimiter=(,right delimiter={)}] {
				\text{\scriptsize  Normal} \& \text{\scriptsize moments} \& \&  \&  \\
				\text{\scriptsize applied to} \& f \&  \& \&  \\
				\&  \&  \&  \& \\
				\&  \&  \&  \& \\
				\&  \&  \& \text{\scriptsize Internal}   \& \text{ \scriptsize moments} \\
				\&  \&  \& \text{\scriptsize applied to} \& g \\
			};
			\end{scope}

			\begin{scope}[shift={(10,0)}]
			\matrix (u) [matrix of math nodes,left delimiter=(,right delimiter={)}] {
				\\
				f\\
				\phantom{f}\\
				\\
				g\\
				\\
			};
			\end{scope}

			\draw (1.6,0) node  [align=left] {$ = $};
			\draw (7,0.5) node  [align=left] {$ 0 $};
			\draw (4,-0.7) node  [align=left] {$ 0 $};
			\draw[decorate,transform canvas={xshift=-1.5em},thick] (b-3-1.south west) -- node[left=2pt] {$\substack{\text{Normal Dofs} \\ \text{values}}$} (b-1-1.north west);
			\draw[decorate,transform canvas={xshift=-1.5em},thick] (b-6-1.south west) -- node[left=2pt] {$\substack{\text{Internal Dofs} \\ \text{values}}$} (b-4-1.north west);
			\end{tikzpicture}}
	\end{center}
	Clearly, there is no interconnection between the function's characterisation on the boundaries and the one performed within the element.
	Thus, showing the {proposition \ref{Prop:NgonsUnisolvence3}} reduces to show independently that
	\IIa{} = 0  \text{ or } \IIb{} = 0\text{ implies } $f|_{\partial K} = 0$
	\text{and}
	$\int_{K}{}{g \, \cdot \, p_k}{x} = 0$, \quad \text{ for all } $p_k\in\mathcal{P}_k$
	 \text{ implies } $g|_{\mathring{K}} = 0.$

	\medskip
\noindent\textit{Step 2.} Let us now consider the boundary characterisation.
	There, by definition of the spaces $\mathcal{H}_{l_1}$ and $\mathcal{H}_{l_2}$, the function $f|_{\partial K}$ is discontinuous at the polytope's vertices and can be decomposed into $n$ vectorial polynomial functions $\{f_j\}_j$ with distinct supports, each of them matching one particular face of the polytope. Thus, we can write
	\begin{equation*}
	f|_{\partial K} = \sum\limits_{j=1}^n r_j\mathbbm{1}_{f_j}
	\end{equation*}
with $r_j\in\bigtimes_{i=1}^d\mathbb{Q}_{\max\{l_1,\,l_2+1\}}(f_j)$ and $f_j$ any face belonging to $\partial K$. With a similar argument than in the previous point, the characterisation of $f|_{\partial K}$ can therefore be done edge-wise, and the determination matrix becomes block - diagonal. We discuss here the characterisation on one particular edge $f_j$ by showing the invertibility of the corresponding matrix block. The arguments naturally transpose to the other ones.

In this perspective, let us show that for any $r_j\in\bigtimes_{i=1}^d\mathbb{Q}_{\max\{l_1,\,l_2+1\}}(f_j),$ it holds
$\{(\IIa)|_{f_j} = 0  \text{ or } (\IIb)|_{f_j} = 0\} \Rightarrow r_j = 0,$
where $(\cdot)|_{f_j}$ represents the subset of the \Dofs  $(\cdot)$ whose support (or evaluation point for point-values) matches (or lies on) $f_j$.

First of all, we recall that on the face $f_j$ the function $r_j$ is a multi-valued polynomial of the form
\begin{equation*}
\resizebox{\textwidth}{!}{	$r_j|_{f_j} = \begin{pmatrix}a_{0,\,1}\\\vdots\\a_{0,\,d}\end{pmatrix}
		+ \sum\limits_{i=\dim(\mathcal{H}_{0})}^{\dim(\mathcal{H}_{l_1}\cap\mathcal{H}_{l_2})}\begin{pmatrix}
		b_{\xi_1{(i)}}+a_{i,\,1}\\\vdots\\b_{\xi_d{(i)}}+a_{i,\,d}						  \end{pmatrix}m_{\alpha_i(x)}
		+ \sum\limits_{i=\dim(\mathcal{H}_{l_1}\cap\mathcal{H}_{l_2})}^{\dim(\mathcal{H}_{l_2})}\begin{pmatrix}
		x_1b_i\\\vdots\\x_db_i
		\end{pmatrix}m_{\alpha_i(x)},$}
\end{equation*}
where $m_{\alpha_i}$ represents a monomial of $\mathbb{Q}_{\max\{l_1,\,l_2\}}$ of multi-index degree $\alpha_i$ such that the set $\{m_{\alpha_i}\}_{i=\dim(\mathcal{H}_{l})}^{\dim(\mathcal{H}_m)}$ forms a base of $\mathcal{H}_{m}\setminus\mathcal{H}_{l}$. Note that the coefficients $\{a_{ij}\}_{i,\,j}$ are defined coordinate-wise while the coefficients $\{b_i\}_i$ are identical for all the components.
	The function $r_j$ is therefore determined by \[\dim(\{\{a_{i,\,m}\}_{\substack{i\in\irange{0}{\dim(\mathcal{H}_{l_1}\setminus\mathcal{H}_0)}\\m\in\irange{1}{d}}},\,\{b_{i}\}_{\substack{i\in\llbracket\dim(\mathcal{H}_{l_1}\setminus\mathcal{H}_0),\\\quad\quad\quad\dim(\mathcal{H}_{l_2})\rrbracket}}\})\]
coefficients.

	As in all configurations the function $r_j$ is determined only through its normal components, let us use the above expression to derive them more specifically. With the normal $n_j = (n_{jx_1,\,\cdots,\,n_{jx_d}})$ to the face $f_j$, it comes
	\begin{equation*}
	\resizebox{\textwidth}{!}{$
		r_j\cdot n_j|_{f_j} = \sum\limits_{m=1}^{d}a_{0,\,m}n_{jx_m}
		+\hspace{-1.5em} \sum\limits_{i=\dim(\mathcal{H}_{0})}^{\dim(\mathcal{H}_{l_1}\cap\mathcal{H}_{l_2})}
		\hspace{-0.2em}\sum\limits_{m=1}^{d}(b_{\xi_m(i)}+a_{i,\,m})n_{jx_m}m_{\alpha_i}(x)
		+\hspace{-1.5em} \sum\limits_{i=\dim(\mathcal{H}_{l_1}\cap\mathcal{H}_{l_2})}^{\dim(\mathcal{H}_{l_2})}\hspace{-1.7em}c_j\,b_i\,m_{\alpha_i}(x),$}
	\end{equation*}
	where $c_j = x\cdot n_j$ is a constant term on the face $f_j$. Reordering the terms, we end up with the formulation
	\begin{equation}\resizebox{0.8\textwidth}{!}{$
		\begin{aligned}
		\label{eq:222}
		r_j\cdot n|_{f_j}& =\sum\limits_{m=1}^{d}\left(\left(a_{0,\,m}+
		\hspace{-1.5em}\sum\limits_{i=\dim(\mathcal{H}_{0})}^{\dim(\mathcal{H}_{l_1}\cap\mathcal{H}_{l_2})}\hspace{-1.2em}a_{i,\,m}m_{\alpha_i(x)}\right)n_{jx_m}\right) \\
		&\qquad+\left(\sum\limits_{m=1}^{d}n_{jx_m}\right) \sum\limits_{i=1}^{\dim(\mathcal{H}_{l_1}\cap\mathcal{H}_{l_2})}\hspace{-1.5em}b_{\xi_m(i)}m_{\alpha_i}(x) + c_j \left(\hspace{-1.5em} \sum\limits_{\hspace{1.5em}i=\dim(\mathcal{H}_{l_1}\cap\mathcal{H}_{l_2})}^{\dim(\mathcal{H}_{l_2})} \hspace{-1.5em}\,b_i\,m_{\alpha_i}(x)\right).
		\end{aligned}$}\end{equation}
	The structure of the retrieved form makes clearly emerge the coefficients that should be used depending on the coordinate-wise behaviour of the polygon.

	In addition, as all the coefficients determining $r_j$ appear in this expression, using \Dofs defined only from the normal components of tested functions is admissible. Thus, the two configurations fitting this framework, we only have to make sure that the set of extracted \Dofs are uniquely characterising each of the involved coefficients. To this aim, we explicit all the possible \Dofs when applied to the function $r_j$. For the sake of clarity, we denote by $\{\sigma_{M_{i,l}}\}_{il}$ the moments designed coordinate-wise, being of the form
	\begin{align*}
	\sigma_{M_{i,l}}&\colon q\mapsto q_{x_i}(x_{jl})n_{jx_i} \quad\quad
	\text{or}\quad\quad \quad\sigma_{M_{i,l}}\colon q\mapsto \SIntL{f_j}{}{q_{x_i}n_{jx_i}p_l}\quad\quad\quad
	\end{align*}
	for some scalar polynomial $p_l$, and by $\{\sigma_{T_l}\}_l$ the ones acting globally, reading
	\begin{align*}
	\sigma_{T_l}&\colon q\mapsto q(x_l)\cdot\,n \quad\quad
	\text{or}\quad\quad \quad	\sigma_{T_l}\colon q\mapsto \SIntL{f_j}{}{q\cdot n\,p_l}\quad\quad
	\end{align*}
	for some scalar polynomial $p_l$. Further, for convenience we denote by $\{\sigma_{V_l}\}_l$ the global \Dofs that comes into play to determining the coordinate - wise coefficients, whose expressions are done in the same way as $\{\sigma_{T_l}\}_l$.
	We now express those \Dofs depending on the coefficients $\{b_{i,\,m}\}_{im}$ and $\{a_{i}\}_{i}$.
	Using the linearity of the \Dofss, plugging the expression \eqref{eq:222} in place of $q$ and setting the permutation operator directly on the multi - indices $\alpha_i$ instead of the coefficients $a_i$, we can rewrite the moments as follows.
	\begin{equation*}
	\begin{split}
		\sigma_{M_{m,l}}\colon (\{a_{i,\,m}\},\,\{b_{i}\})\longmapsto &
	\, b_{0,\,m}\int_{f_j}{}{n_{jx_m}\,p_l}	+\hspace{-1.5em}\sum\limits_{i=\mathcal{H}_{0}}^{\dim(\mathcal{H}_{l_1}\cap\mathcal{H}_{l_2})}\hspace{-1.2em}a_{i,\,m}\int_{f_j}{}{m_{\alpha_i}(x)n_{jx_m}\,p_l} \\
	+&\hspace{-1.5em} \sum\limits_{i=\mathcal{H}_{0}}^{\dim(\mathcal{H}_{l_1}\cap\mathcal{H}_{l_2})}\hspace{-1.2em}b_i\,
	\int_{f_j}{}{n_{jx_m}m_{\xi_m{(\alpha_i)}}(x)p_l}
	+\hspace{-1.5em} \sum\limits_{i=\dim(\mathcal{H}_{l_1}\cap\mathcal{H}_{l_2})}^{\dim(\mathcal{H}_{l_2})}\hspace{-1.2em}b_i\,
	\int_{f_j}{}{(x_m\,n_{jx_m}m_{\alpha_i}(x))p_l}
	\end{split}
	\end{equation*}
	\begin{equation*}
	\begin{split}
		\sigma_{T_l}\colon (\{a_{i,\,m}\},\,\{b_{i}\})\longmapsto &
	\sum\limits_{m=1}^{d}a_{0,\,m}\int_{f_j}{}{n_{jx_m}p_l}
	+\sum\limits_{m=1}^{d}\hspace{-0.5em}\sum\limits_{i=\dim(\mathcal{H}_{0})}^{\dim(\mathcal{H}_{l_1}\cap\mathcal{H}_{l_2})}\hspace{-1.2em}a_{i,\,m}\int_{f_j}{}{m_{\alpha_i}(x)n_{jx_m}p_l}
 \\
	+&\hspace{-1.2em}\sum\limits_{i=\mathcal{H}_{0}}^{\dim(\mathcal{H}_{l_1}\cap\mathcal{H}_{l_2})}
	\hspace{-1.2em}  b_i
	\sum\limits_{m=1}^{d}
	 \left(\int_{f_j}n_{jx_m}m_{\xi_m(\alpha_i)}(x)p_l\right)
	+ \hspace{-1.2em}\sum\limits_{i=\dim(\mathcal{H}_{l_1}\cap\mathcal{H}_{l_2})}^{\dim(\mathcal{H}_{l_2})}\hspace{-1.2em}b_i\,
	\int_{f_j}{}{\big(c_jm_{\alpha_i}(x)p_l\big)}
	\end{split}
	\end{equation*}
	Thus, defining the component-wise parts of the global moments $\sigma_{T_{l}}$  by $\sigma_{T_{m,l}}(q) = \int_{f_j}{n_{jx_m}q\,p_l}$ such that $\sigma_{T_{l}} = \sum_{m=1}^d\sigma_{T_{m,l}}$, one can express any \Dofs of the two considered sets as
	\begin{align*}
	\sigma_{M_{m,l}}\colon q\longmapsto
	\,& a_{0,\,m}\sigma_{M_{m,l}}(1) + \hspace{-1.5em}\sum\limits_{i=\mathcal{H}_{0}}^{\dim(\mathcal{H}_{l_1}\cap\mathcal{H}_{l_2})}\hspace{-1.2em}a_{i,\,m}\sigma_{M_{m,l}}(m_{\alpha_i})\\
	&\hspace{0.2em}+\hspace{-1.5em} \sum\limits_{i=\mathcal{H}_{0}}^{\dim(\mathcal{H}_{l_1}\cap\mathcal{H}_{l_2})}\hspace{-1.2em}b_i\,
	\sigma_{M_{m,l}}(m_{\xi_m{(\alpha_i)}})+\hspace{-1.5em} \sum\limits_{i=\dim(\mathcal{H}_{l_1}\cap\mathcal{H}_{l_2})}^{\dim(\mathcal{H}_{l_2})}\hspace{-1.2em}b_i\,
	\sigma_{M_{m,l}}(x_m m_{\alpha_i})
	\end{align*}
	\begin{align*}
	\text{and \quad}\sigma_{T_l}\colon q\longmapsto &
	\sum\limits_{m=1}^{d}a_{0,\,m}\sigma_{T_{m,l}}(1)+\hspace{-1.5em}\sum\limits_{i=\mathcal{H}_{0}}^{\dim(\mathcal{H}_{l_1}\cap\mathcal{H}_{l_2})}\hspace{-0.6em}\sum\limits_{m=1}^{d}a_{i,\,m}\sigma_{T_{m,l}}(m_{\alpha_i}) \\
	&\hspace{0.2em}+\hspace{-1.5em}\sum\limits_{i=\mathcal{H}_{0}}^{\dim(\mathcal{H}_{l_1}\cap\mathcal{H}_{l_2})}\hspace{-1.2em}  b_i\sum\limits_{m=1}^{d} \sigma_{T_{m,l}}(m_{\xi_m({\alpha_i})})
	+ \hspace{-1.2em}\sum\limits_{i=\dim(\mathcal{H}_{l_1}\cap\mathcal{H}_{l_2})}^{\dim(\mathcal{H}_{l_2})}\hspace{-1.2em}a_ic_j\,\sigma_{T_{l}}(m_{\alpha_i}).
	\end{align*}
	Note that in view of deriving the determination matrix, the last term can also be decomposed as follows. \[\sum\limits_{i=\dim(\mathcal{H}_{l_1}\cap\mathcal{H}_{l_2})}^{\dim(\mathcal{H}_{l_2})}\hspace{-1.2em}b_ic_j\,\sigma_{T_{l}}(m_{\alpha_i}) = \sum\limits_{i=\dim(\mathcal{H}_{l_1}\cap\mathcal{H}_{l_2})}^{\dim(\mathcal{H}_{l_2})}\hspace{-1.2em}b_i\,\sum\limits_{m=1}^{d}\sigma_{T_{l,m}}(x_mm_{\alpha_i}).\]
	Similar relations for $\sigma_{V}$ can ve derived from the expression of $\sigma_{T}$. Thus, we can rewrite the \Dofs as a dot product and derive the characterisation matrix $\Sigma$
	\[
	(\sigma_{M_{1,1}},\,\cdots ,\,\sigma_{M_{d,l}},\,
	\sigma_{V_{1}},\,\cdots ,\,\sigma_{V_{l}},\,
	\sigma_{T_{1}},\,\cdots,\,\sigma_{T_{l}},\,)^T
	= \Sigma(\{a_{i,\,m}\}_{im},\,\{b_{i}\}_i)^T
	\]
	which shape is given in the figure \ref{enormeMatrice}. We now investigate its structure.

	First of all, as the number of extracted \Dofs from the two sets \IIa{} and \IIb{} matches the number of coefficients determining $r_j$, the matrix $\Sigma$ is a square matrix.

	Let us focus on the top two-by-two left blocks, surrounded in blue. They correspond to the coefficients that should be determined coordinate-wise.
	Thus, by construction, there are $\dim(\{a_{im}\}_{im}) = d(l_1+1)^{d-1}$ columns. And by the definition of the configurations \I{} and \II{}, we have $\dim(\{\sigma_{M_{i,\,j}}\}_{ij}) = d + d((l_1+1)^{d-1}-1) = d(l_1+1)^{d-1}$. Therefore, this submatrix is a square matrix.
	Furthermore, each subblock corresponds to one member of the decomposition of $q$ tested through coordinate - wise \Dofs whose kernels are built on the same monomial. Therefore, as the \Dofs $\{\sigma_{M_{i,\,j}} \}$ consider one normal component only, the coefficients $\{a_{im}\}_{im}$ for $m\neq j$ are not involved, and the subblocks are diagonal. Thus, those submatrices are invertible and in particular their columns and rows are linearly independent.

	On the other side of the matrix, the last bottom block surrounded in deep red matches the Raviart-Thomas moments tuning members of $\Hk|_{\partial K}$ living exclusively in $x B_k|_{f_j}$. It is then a submatrix of the classical Raviart-Thomas' one, and is thus invertible. In particular, its rows and columns are linearly independent.

	The extended bottom right submatrix highlighted in dashed red corresponds to the previously described high-order submatrix of the Raviart-Thomas's setting, enriched by the moments $\{\sigma_V\}$ tuning the behaviour of members of $\Hk|_{f_j}$ falling in the intersection $\mathcal{H}_{l_1}\cap\mathcal{H}_{l_2}(f_j)$.

	This matrix is equivalent to the full Raviart-Thomas setting. Indeed, even if the moments $\{\sigma_V \}$ have to be slightly modified from the Raviart-Thomas setting in the configuration $I$, this modification leaves the projection order unchanged and the integrand still belongs to $\Hk|_{f_j}$. Therefore, the dashed line block is invertible and its columns and rows are linearly independent.

	There is only left to show that there is no linear dependence between rows of different row blocks. As the \Dofs are linear forms, it is enough to show that the integrand of moments (or  polynomials constructing the point - wise values) that involve the same monomial are linearly independent.

	Indeed, being linear forms whose kernels are polynomials, the \Dofs can combine each other only if their integrand ($q$ tested against the kernel) involve -- up to constants -- the same monomials. We then have to show that in both configurations, the rows involving terms whose projection onto the kernel can be expressed from a same monomial are linearly independent.

	In the configuration \II, this property comes automatically. Indeed, the only interaction between \Dofs having integrands sharing the same monomial order (and then possibly being based on the same monomial) is possible between \eqref{eq:Doffs3} and \eqref{eq:Doffs2} when $|p_k|=l_1+1$. Indeed, by definition of $\Hk$, the polynomial $p_k\cdot n$ in \eqref{eq:Doffs3} is only of order $l_1$. However, no combination of \eqref{eq:Doffs2} can form the moments \eqref{eq:Doffs3}. Indeed, for any real coefficients $c_i$ it holds
	\[
	\sum c_i q_{x_i}n_{x_i}x_i^{l_1+1}\not\equiv q\cdot n p_k
	\]
	for any monomial $p_k$ such that $|p_k| = l_1+1$.
	Note that in the left hand side, all the $c_i$ should be non-null to reconstruct the term $p\cdot n$. However, doing so no factorisation by a single monomial such that
	\[
	\sum_{i=1}^{d} c_i q_{x_i}n_{x_i}x_i^{l_1+1} = \left(\sum_{i=1}^{d} c_i q_{x_i}n_{x_i}\right) p_k
	\] is possible.

	\begin{figure}[h]
		\begin{center}
		\includegraphics[height=\textheight]{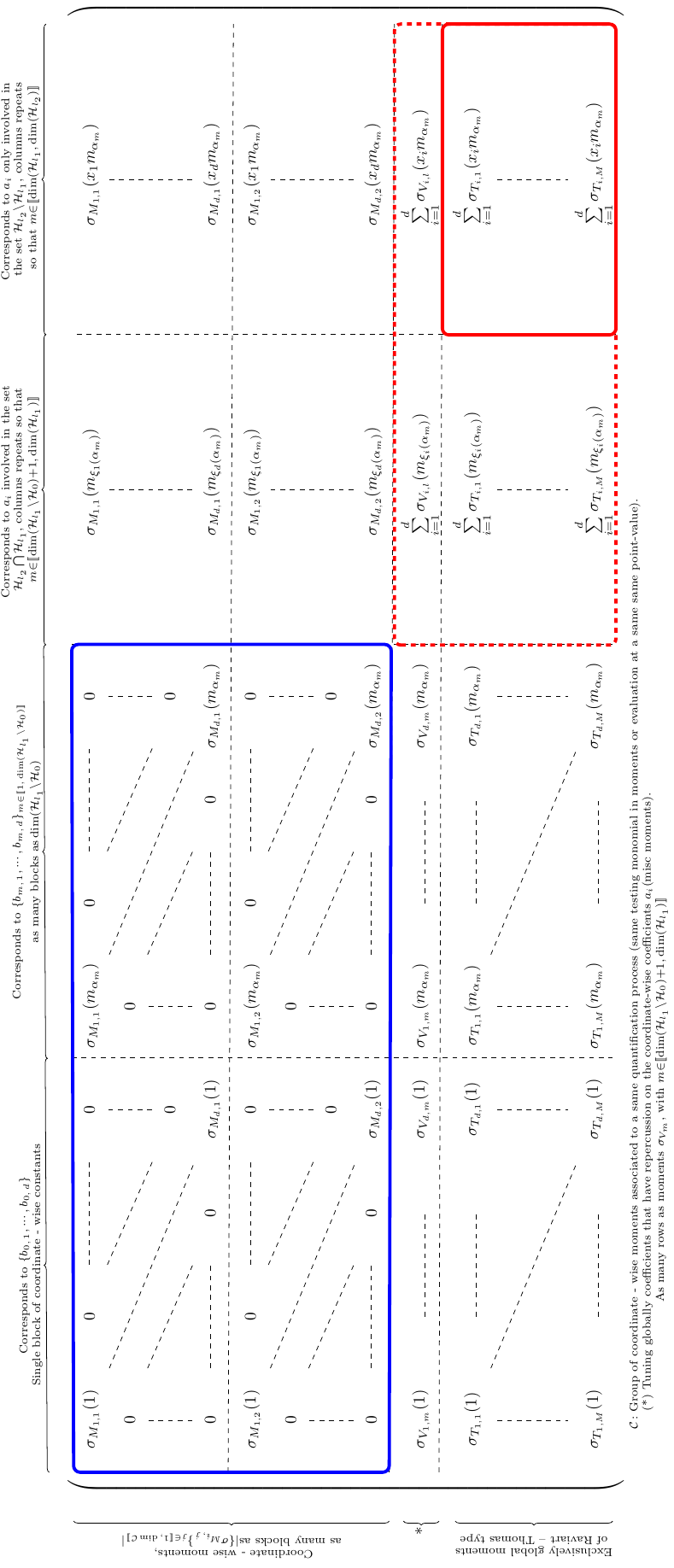}
		\end{center}
		\caption{\label{enormeMatrice} Disposition of the degrees of freedom.}
	\end{figure}

	Thus, the designed moments are linearly independent, and no row combination can occur for any tested polynomial belonging to $\Hk|_{\partial K}$.


	All in all, for both configurations all the rows are linearly independent. As by construction we have as many relations as unknowns, the matrix is invertible. Thus, we get a null kernel, meaningly
	\[\IIa = 0  \text{ or } \IIb = 0 \text{ implies } f|_{\partial K} = 0.\]

	\noindent\textit{Step 3.} Let us now consider the internal characterisation of functions living in $\Hk$. From the first point of the proof, it is enough to study the characterisation of $g$ within the inner cell. By definition of $\Hk$, any function $g\in\Hk|_{\mathring{K}}$ can be decomposed over a set of Poisson's solutions as follows.
	\[
	g = \sum\limits_{i=1}^d \sum\limits_{i=1}^{\dim A_k} a_{i,\,j}e_ju_i + \sum\limits_{i=1}^{\dim B_k}b_j x \tilde{u}_i
	\]
	Here, the vector $e_j$ stands for $e_j = (0,\,\cdots,\,1,\,0,\,\cdots,\,0)^T$ where the $1$ is in the $j^{\mathrm{th}}$ position. The functions $u_i$ and $\tilde{u}_i$ represents the Poisson's solutions of the problems
	\begin{equation}\label{eq:form}
	\begin{cases}
	\Delta u_i=p_i,\,p_i\in\mathbb{Q}_{m_1}(K)\\
	u_i|_{\partial K} = 0
	\end{cases}\,\,\text{and}\quad
	\begin{cases}
	\Delta \tilde{u_i}=\tilde{p_i},\,\tilde{p_i}\in\mathbb{Q}_{[m_2]}(K)\\
	\tilde{u_i}|_{\partial K} = 0
	\end{cases}
	\end{equation}
	where $\{p_i\}_i$ and $\{\tilde{p_i}\}_i$ form respectively a basis of $\mathbb{Q}_{m_1}(K)$  and $\mathbb{Q}_{[m_2]}(K)$.
	In any presented definition of the \Dofss, the internal characterisation is done through moment - based \Dofs of the form
	\[
	\sigma_{I_k}\colon	q\mapsto \, \Int{K}{}{q\cdot\, p_k}{x}
	\]
	where the kernels $p_k\in\mathcal{P}_k$ consist of linearly independent polynomials belonging to $\mathbb{Q}_{\max\{m_1,\,m_2+1\}}(K)$, or of the solution of their corresponding problems of the form \eqref{eq:form}. Therefore, we can derive a characterisation matrix in the same spirit as in the case of the normal characterisation.
	\begin{equation}\label{matrmat}
		\hspace{10em}\resizebox{0.70\textwidth}{!}{
			\SmallMatrix{1.4cm}{1.4cm}}
	\end{equation}

	Let us consider the case where $\mathcal{P}_k$ forms a polynomial projection space. There, none of the $p_k\in\mathcal{P}_k$ is the zero function. In the same time, the functions $\{\{u_i\}_i,\,\{x\tilde{u_i}\}_i\}$ are linearly independent, and being solutions to some Poisson's problem with non-zero second member, they are by construction not identically vanishing on $K$. Indeed, even when $m_2<m_1$ where second members of the problems \eqref{eq:form} lives both in $\mathbb{Q}_{m_1}$ and $\mathbb{Q}_{[m_2]}$, it holds $\Delta (x\,\tilde{u_i}) = 2\nabla \cdot \tilde{u_i} + \Delta(\tilde{u_i})$. Thus, it is impossible to combine linearly the function $x\,\tilde{u_i}$ with functions of the set $\{u_i\}_i$.

	Furthermore, the degrees of the polynomials belonging to the space $\mathcal{P}_k$ are lower or equal than the highest degree of the second members of the Poisson's problem defining the space $\Hk$. Thus, every projection of function of $\Hk$ onto the space $\mathcal{P}_k$ is not null. And as the internal moments are linear forms, any linear combination of those moments at fixed $p_k$ could have its integrand factorised by the kernel $p_k$ for any $p_k\in\mathcal{P}_k$, transferring the linear independency of the set $\{\{u_i\}_i,\,\{x\tilde{u_i}\}_i\}$ to the terms $\{\int u_i\cdot p_k\}_{i}$ for any fixed $p_k\in\mathcal{P}_k$.

	Lastly, as the space $\mathcal{P}_k$ contains only linearly independent functions the previous argument can be repeated for each row of the matrix defined in \eqref{matrmat}. And as by construction the number of internal \Dofs matches the dimension of the space $\Hk|_{\mathring{K}}$, the linear independence of functions of $\mathcal{P}_k$ combined with the linear independence of the tested functions transfer automatically to the moments tested against a basis of $\Hk|_{\mathring{K}}$. Thus, the internal submatrix is invertible. The same reasoning can be applied when $\mathcal{P}_k$ is built from the Poisson's solutions themselves, as the projections of functions would decompose the functions directly.

	Merging the above points together, we get that
	$\IIa = 0  \text{ or } \IIb = 0$  implies $f|_{\partial K} = 0$ and
	$\int_{K}{}{g \, \cdot \, p_k}{x} = 0$ { for all } $p_k\in\mathcal{P}_k$ implies $g|_{\mathring{K}} = 0$.
	From this, we get that  for $q\in\Hk$
	$\IIa = 0  \text{ or } \IIb = 0 $ and
	$\int_{K}{}{q \, \cdot \, p_k}{x} = 0$ { for all } $p_k\in\mathcal{P}_k$
	implies $q=0$.

\end{proof}

\begin{proof}[Proof of the Propositions \ref{GeneralProof2} and \ref{Prop:NgonsUnisolvence3}]
	\label{theproof} The proof of the {proposition \ref{GeneralProof2}} is a strai\-ght\-forward ge\-ne\-ra\-lisation of the one presented for the two examples \IIa{} and \IIb{} in the lemma \ref{LemmaImportant}. The only change lies in the extraction of the \Dofss, which impacts the matrix only on the top left two by two blocks describing the coordinate-wise behaviours.
	As the extraction fulfils the assumptions \ref{Ass:Dofs0}, the rows involving terms whose projection onto the kernel can be expressed from a same monomial are linearly independent. Thus, the same arguments as above can be applied and the conclusion follows.

	In particular, by this admissibility criterion there cannot be more than $d+1$ polynomials reducing to the same moments' kernel or to an equivalent point-value quantifier. Thus, as we have $d+1$ coordinate-wise moments to tune per decomposed monomial, there is no over-determination at a fixed polynomial degree. The constraint on the extraction of \Dofs ensures the non over-determination overall. Further, the linear independence of the sub-matrix's columns is ensured as those polynomials cannot be linearly dependent. Thus, by linearity of the \Dofss, the independence of the kernels transfers to the moments and there is no row dependency. The submatrix block corresponding to any specific order is therefore invertible, and the same conclusion as in the {proof \ref{LemmaImportant}} follows.

	The proposition \ref{GeneralProof2} thus holds by the number of \Dofss, matching the dimension of the space $\Hk$. Indeed, as by the proposition \ref{Prop:NgonsUnisolvence3} the kernel of the linear operator defined by the set of \Dofs has a null kernel providing their unisolvence when enclosed within the space $\Hk$.
\end{proof}

\begin{proof}[Proof {of the relation \eqref{Prop:DimNgons}}]
	As for any $l_1\leq 0$ the two natural subspaces are in direct sum, recalling the block construction of $\Hk$ allows the dimension of the space $\Hk$ to be easily derived. We can simply add the dimension of the two main subspaces $(A_k)^d$ and $x\,B_k$ to retrieve the dimension of $\Hk$. Let us derive their respective dimensions.

First, we compute the dimension of $A_k$. In the way  presented in \cite{VEIGA2013}, one can get it by using the superposition theorem. Indeed, for any second member belonging to $\mathbb{Q}_{m_1}$ and any boundary function $p_k\mathbbm{1}_{f}\in L^2(K)$, there exists a unique solution to the Poisson's problems defining $A_k$ ({see e.g. \cite{chabrowski2006dirichlet}}). Thus, reading out the structure of the set $A_k$ implies the following relation.
	\begin{align*}
	\dim A_k &= \dim \mathcal{H}_{l_1}(\partial K) + \dim \mathbb{Q}_{m_1}(K)\\
	&= n(l_1+1)^{d-1} + (m_1+1)^d
	\end{align*}
Therefore, as $(A_k)^d$ is a simple Cartesian product of $d$ copies of $A_k$, we have immediately $\dim A_k = d(\dim A_k) = d(n(l_1+1)^{d-1} + (m_1+1)^d)$. In the exact same way, we retrieve the dimension of $B_k$ by
	\begin{align*}
	\dim B_k &= \dim \mathcal{H}_{l_2}(\partial K) + \dim \mathbb{Q}_{[m_2]}(K)\\
	&= n(l_2+1)^{d-1} + (m_2+1)^d - m_2^d.
	\end{align*}
Last, we recall that the space $x\, B_k$ simply corresponds to an identical $d$ - duplication of the space $B_k$ where each coordinate has been multiplied by the corresponding spatial variable. Thus, there is no liberty adjunction during its construction, and the dimension of $x\, B_k$ equals the one of $B_k$. By combining this, we finally get
	\begin{align*}
	\dim \Hk &= d\dim A_k + \dim B_k\\
	&= d(n(l_1+1)^{d-1} + (m_1+1)^d) + n(l_2+1)^{d-1} + (m_2+1)^d - m_2^d.
	\end{align*}
\end{proof}

Finally, we point out the relation between the proposed spaces and the slight restriction of the one introduced in \cite{da2016h}. Indeed, we focus on relation \eqref{eq:connection}.

\begin{proof}[Sketch of the proof that $\Edit{\tilde{\mathcal{V}}}_{2,\,k}(K)\subset \mathbb{H}_k(K)$]
Let us quickly show that any element of $\Edit{\tilde{\mathcal{V}}}_{2,\,k}(K)$ can be recast as an element of $\mathbb{H}_k(K)$ for the coefficients $(l_1,\,l_2) = (0,\,k)$ and $(m_1,\,m_2) = (k-1,\,-1)$.
\vspace{0.5em}

\noindent\textit{Regularity.}
Any element $v\in\Edit{\tilde{\mathcal{V}}}_{2,\,k}(K)$ belongs to $H^1(K)$. Indeed, $v\in H(\mathrm{div},K)\cap H(\mathrm{rot},K)\subset (H^1(K))^d_.$ Furthermore, as any element $w$ of $\mathbb{H}_k(K)$ writes $w = u_1 + x\,u_2$ with $u_1,\,u_2$ both having for regularity $ H^1(K)$, $w\in H^1_{loc}(K)$. And since $K$ is compact and bounded, $H^1_{loc}(K)=H^1(K)$. Thus, the regularity asked for any element $v$ to be in $\Edit{\tilde{\mathcal{V}}}_{2,\,k}(K)$ is stricter than the one asked for any element $w$ to be in $\mathbb{H}_k(K)$.
\vspace{0.5em}

\noindent\textit{On the boundary.}
Any $v\in\Edit{\tilde{\mathcal{V}}_{2,\,k}(K)}$ \Edit{is a polynomial} and satisfies $v\cdot n|_f \in \mathbb{P}_f(K)$ on every face $K$ of the polygon $K$. But $v\cdot n|_f $ is nothing else than the linear combination
$\sum_{i=1}^d v_in_{x_i}$ of the coordinate-wise functions $v_i$ with the normal's coefficients $n_{x_i}$. Thus, each polynomial $v_i$ has no choice but to live in the space $\mathbb{P}_{k,\,\cdots,\,k,\,k+1,\,k,\cdots,\,k}(K)$, where the $k+1$ is in the $i^\mathrm{th}$ position. The space $\times_{i=1}^d\mathbb{P}_{k,\,\cdots,\,k,\,k+1,\,k,\cdots,\,k}(K)$ is indeed the smallest space that is polynomial (required by the linear combination) and that contains all the functions $v$ such that $v\cdot n|_f \in \mathbb{P}_f(K)$. Note that allowing a higher degree in the $i^\mathrm{th}$ variable is required as on each face $f\in\partial K$, $x\cdot n \equiv c$ for some constant $c$. And $\times_{i=1}^d\mathbb{P}_{k,\,\cdots,\,k,\,k+1,\,k,\cdots,\,k}(K) \subset \mathbb{P}_0(K)+x\,\mathbb{P}_{k}(K)$, which is exactly the structure of the $\mathbb{H}_k(K)$ space on the boundary.
\vspace{0.5em}

\noindent\textit{Within the element.}
For any  $v\in\Edit{\tilde{\mathcal{V}}}_{2,\,k}(K)$, it holds
\[
\begin{cases}
\nabla (\nabla \cdot v) &\in \nabla (\mathbb{P}_{k-1}(K))\\
\nabla \times v &\in  \mathbb{P}_{k-1}(K)
\end{cases}
\]
Thus, it comes $\underbrace{\nabla(\nabla \cdot v)}_{\in\nabla(\mathbb{P}_{k-1}(K))}-\underbrace{\nabla\times(\nabla \times v)}_{\in\nabla(\mathbb{P}_{k-\Edit{2}}(K))}\in \nabla(\mathbb{P}_{k-1}(K))$, which implies $\nabla^2v\in \nabla(\mathbb{P}_{k-1}(K))$ and writes
\[
\nabla^2v = \begin{pmatrix}
\sum_{i=1}^{d}  \frac{\partial^2v_1}{\partial x_i^2}\\
\vdots\\
\sum_{i=1}^{d}  \frac{\partial^2v_d}{\partial x_i^2}
\end{pmatrix}
\in
\begin{pmatrix}
\mathbb{P}_{k-2,\,k-1,\,\cdots,\,k-1}(K)\\
\vdots\\
\mathbb{P}_{k-1,\,k-1,\,\cdots,\,k-2}(K)
\end{pmatrix}
\]
Thus, we have naturally
\[
\nabla^2v \in \Big\{ u,\, \begin{pmatrix}\Delta u_1\\\vdots\\ \Delta u_d\end{pmatrix},\, \Delta u_i\in\mathbb{Q}_{k-1}(K)\Big\}.
\]
\vspace{0.5em}

\noindent\textit{Summary}
Let $v\in\mathcal{V}_{2,\,k}(K)$. Then $v$ can be decomposed as follows:
\[
v = (\mathring{v_1}+\mathring{v_2})|_{\mathring K} + (\bar{v_1}+\bar{v_2})|_{\partial K}
\]
where
\[
\begin{cases}
\bar{v_1} \in (\mathbb{P}_0(\partial K))^d \quad &\text{represents the constant part of } v \text{ on the boundary}\\
\bar{v_2} \in (\mathbb{P}_k(\partial K)\setminus\mathbb{P}_0(\partial K))^d \quad &\text{represents the higher parts of } v \text{ on the boundary}
\end{cases}
\]
and
\[
\begin{cases}
\mathring{v_1} = v|_{\mathring{K}} - p_{[1],\,\cdots,\,[1]} \\
\bar{v_2} = p_{[1],\,\cdots,\,[1]}
\end{cases}
\]
for any $ p_{[1],\,\cdots,\,[1]} $ belonging to $x\mathbb{P}_{0}$. Doing so, we get
\[
\begin{cases}
\bar{v_1}|_{\partial K} & \in (\mathbb{P}_0(\partial K))^d\\
\mathring{v_1}|_{\mathring{K}} & = v|_{\mathring{K}} - p_{[1],\,\cdots,\,[1]}
\end{cases}
\]
and
\[
\begin{cases}
\bar{v_2}|_{\partial K} &  (\mathbb{P}_k(\partial K)\setminus\mathbb{P}_0(\partial K))^d\\
\mathring{v_2}|_{\mathring{K}} & = p_{[1],\,\cdots,\,[1]}
\end{cases}
\]
There, we get straightforwardly from the previous paragraphs:
\[
\begin{cases}
v_1 \in H^1(K)\\
v_1\cdot n|_{\partial K} \in \mathbb{P}_0(\partial K) \subset \mathbb{Q}_0(\partial K)\\
\Delta v_1|_{\mathring{K}} = \underbrace{v|_{\mathring{K}}}_{\substack{\in\nabla(\mathbb{P}_{k-1}(K))\\\subset\mathbb{P}_{k-1}(K)}}-\underbrace{\Delta(p)}_{=0} \quad \in \mathbb{P}_{k-1}(K)\subset\mathbb{Q}_{k-1}(K)
\end{cases}
\]
and
\[
\begin{cases}
v_2 \in H^1(K)\\
v_2\cdot n|_{\partial K} = v|_{\partial K}\cdot n - v_1\cdot n|_{\partial K} \quad &\in \mathbb{P}_k(\partial K) \subset \mathbb{Q}_k(\partial K)\\
\Delta v_2|_{\mathring{K}} = \Delta(x\, c)  = x\Delta(c) = 0&\in \mathbb{P}_{[-1]}(K)\subset \mathbb{Q}_{[-1]}(K)
\end{cases}
\]
for some constant $c\in\mathbb{P}_0(K)$.  Writing $v_2 $ as $v_2 = x w$ on the boundary with $w \in x\,\mathbb{P}_{k,\,\cdots,\,k}(\partial K)$ \newline$ \subset \times_{i=1}^d\mathbb{P}_{k,\,\cdots,\,k,\,k+1,\,k,\cdots,\,k}(\partial K) \setminus \mathbb{P}_0(\partial K)$, it comes further
\[
\begin{cases}
w \in H^1(K)\\
w\cdot n|_{\partial K} &\in \mathbb{P}_{k}(\partial K) \subset \mathbb{Q}_k(\partial K)\\
\Delta w|_{\mathring{K}} = 0&\in \mathbb{P}_{-1}(K)\subset \mathbb{Q}_{[-1]}(K)
\end{cases}
\]
Taking  without loss of generality $c=1$, considering $v = v_1 + x\,w$ and setting $(l_1,\,l_2) = (0,\,k)$, $(m_1,\,m_2)=(k-1, \, -1)$, we get:
\[ v_1\in \mathcal{A}_k,\quad v_2 = x\,w \in x\,\mathcal{B}_k,\]
and therefore,
\[
v = v_1 + v_2 \in \mathbb{H}_k(K).
\]
%
\end{proof}

	\bibliographystyle{habbrv}
	\bibliography{Refs}

\begin{thebibliography}{10}

\bibitem{abgrall:hal-01820176}
R.~Abgrall, E.~Le~M{\'e}l{\'e}do, and P.~{\"O}ffner.
\newblock On the connection between residual distribution schemes and flux
  reconstruction.
\newblock 2018.

\bibitem{abgrall2019class}
R.~Abgrall, {\'E}.~Le~M{\'e}l{\'e}do, and P.~{\"O}ffner.
\newblock {A class of finite dimensional spaces and {$H$-(div)} conformal
  elements on general polytopes}.
\newblock 2019.

\bibitem{zbMATH06823756}
J.~Aghili, D.~Di~Pietro, and B.~Ruffini.
\newblock An \(hp\)-hybrid high-order method for variable diffusion on general
  meshes.
\newblock {\em Computational Methods in Applied Mathematics}, 17(3):359--376,
  2017.

\bibitem{refId0}
{Beir\~ao da Veiga, Louren\c{c}o}, {Brezzi, Franco}, {Marini, Luisa Donatella},
  and {Russo, Alessandro}.
\newblock Mixed virtual element methods for general second order elliptic
  problems on polygonal meshes.
\newblock {\em ESAIM: M2AN}, 50(3):727--747, 2016.

\bibitem{zbMATH06966728}
F.~Bonaldi, D.~Di~Pietro, G.~Geymonat, and F.~Krasucki.
\newblock A hybrid high-order method for kirchhoff-love plate bending problems.
\newblock {\em ESAIM: Mathematical Modelling and Numerical Analysis},
  52(2):393--421, 2018.

\bibitem{zbMATH07078646}
L.~Botti, D.~Di~Pietro, and J.~Droniou.
\newblock A hybrid high-order method for the incompressible navier-stokes
  equations based on temam's device.
\newblock {\em Journal of Computational Physics}, 376:786--816, 2019.

\bibitem{Brezzi1985}
F.~Brezzi, J.~Douglas, and L.~D. Marini.
\newblock Two families of mixed finite elements for second order elliptic
  problems.
\newblock {\em Numerische Mathematik}, 47(2):217--235, 1985.

\bibitem{chabrowski2006dirichlet}
J.~Chabrowski.
\newblock {\em The Dirichlet problem with L2-boundary data for elliptic linear
  equations}.
\newblock Number 1482. Springer, 2006.

\bibitem{chen2017minimal}
W.~Chen and Y.~Wang.
\newblock {Minimal degree H(curl) and H(div) conforming finite elements on
  polytopal meshes}.
\newblock {\em Mathematics of Computation}, 86(307):2053--2087, 2017.

\bibitem{cockburn2012discontinuous}
B.~Cockburn, G.~E. Karniadakis, and C.~W. Shu.
\newblock {\em Discontinuous Galerkin methods: theory, computation and
  applications}, volume~11.
\newblock Springer Science \& Business Media, 2012.

\bibitem{VEIGA2013}
L.~B. Da~Veiga, F.~Brezzi, A.~Cangiani, G.~Manzini, M.~L. D., and A.~Russo.
\newblock Basic principles of virtual element methods.
\newblock {\em Mathematical Models and Methods in Applied Sciences},
  23(01):199--214, 2013.

\bibitem{da2016h}
L.~B. Da~Veiga, F.~Brezzi, L.~D. Marini, and A.~Russo.
\newblock {H(div) and H (curl)-conforming virtual element methods}.
\newblock {\em Numerische Mathematik}, 133(2):303--332, 2016.

\bibitem{zbMATH06596741}
L.~B. Da~Veiga, F.~Brezzi, L.~D. Marini, and A.~Russo.
\newblock Mixed virtual element methods for general second order elliptic
  problems on polygonal meshes.
\newblock {\em ESAIM: Mathematical Modelling and Numerical Analysis},
  50(3):727--747, 2016.

\bibitem{dipietro2014}
D.~Di~Pietro and S.~Lemaire.
\newblock An extension of the crouzeix-raviart space to general meshes with
  application to quasi-incompressible linear elasticity and stokes flow.
\newblock {\em Mathematics of Computation}, 84(291):1--31, 2015.

\bibitem{10.1093/imanum/drw003}
D.~A. Di~Pietro and A.~Ern.
\newblock {Arbitrary-order mixed methods for heterogeneous anisotropic
  diffusion on general meshes}.
\newblock {\em IMA Journal of Numerical Analysis}, 37(1):40--63, 05 2016.

\bibitem{Dubois2017}
F.~Dubois, I.~Greff, and C.~Pierre.
\newblock Raviart-thomas finite elements of petrov-galerkin type.
\newblock {\em ESAIM. Mathematical Modelling and Numerical Analysis},
  53(5):1553–1576, 2017.

\bibitem{Gillette2014}
A.~Gillette, R.~A., and B.~C.
\newblock Construction of scalar and vector finite element families on
  polygonal and polyhedral meshes.
\newblock {\em Computational Methods in Applied Mathematics}, 16(4):667–683,
  2016.

\bibitem{huynh2007flux}
H.~T. Huynh.
\newblock A flux reconstruction approach to high-order schemes including
  discontinuous galerkin methods.
\newblock In {\em 18th AIAA Computational Fluid Dynamics Conference}, page
  4079, 2007.

\bibitem{lesaint}
P.~Lesaint and P.~Raviart.
\newblock {\em {Mathematical aspects of Finite Element in Partial Differential
  Equations}}, chapter {On a Finite Element Method for Solving the neutron
  transport equation}, pages 89--123.
\newblock Academic Press, 1974.
\newblock C. De Boor, ed.

\bibitem{zbMATH06666913}
K.~Lipnikov, G.~Manzini, and M.~Shashkov.
\newblock Mimetic finite difference method.
\newblock {\em Journal of Computational Physics}, 257:1163--1227, 2014.

\bibitem{zbMATH05897522}
R.~Loub\`ere, P.~H. Maire, and M.~Shashkov.
\newblock Reale: a reconnection arbitrary-lagrangian-eulerian method in
  cylindrical geometry.
\newblock {\em Computers and Fluids}, 46(1):59--69, 2011.

\bibitem{Nedelec1980}
J.~N{\'e}d{\'e}lec.
\newblock Mixed finite elements in {$\mathbb{R}^3$}.
\newblock {\em Numerische Mathematik}, 35(3):315--341, 1980.

\bibitem{ranocha2016summation}
H.~Ranocha, P.~{\"O}ffner, and T.~Sonar.
\newblock Summation-by-parts operators for correction procedure via
  reconstruction.
\newblock {\em Journal of Computational Physics}, 311:299--328, 2016.

\bibitem{RTOriginal}
P.~A. Raviart and J.~M. Thomas.
\newblock A mixed finite element method for 2-nd order elliptic problems.
\newblock In {\em Mathematical Aspects of Finite Element Methods}, pages
  292--315, Berlin, Heidelberg, 1977. Springer Berlin Heidelberg.

\bibitem{reeds}
W.~Reeds and T.~Hill.
\newblock Triangular mesh methods for the neutron transport equation.
\newblock Technical Report LA-UR-73-479, Los Alamos, 1973.

\bibitem{Sukumar2006}
N.~Sukumar and E.~A. Malsch.
\newblock Recent advances in the construction of polygonal finite element
  interpolants.
\newblock {\em Archives of Computational Methods in Engineering}, 13(1):129,
  2006.

\bibitem{Talischi2013}
C.~Talischi, A.~Pereira, G.~H. Paulino, I.~F.~M. Menezes, and M.~S. Carvalho.
\newblock Polygonal finite elements for incompressible fluid flow.
\newblock {\em International Journal for Numerical Methods in Fluids},
  74(2):134--151, 2013.

\bibitem{vincent2011new}
P.~Vincent, P.~Castonguay, and A.~Jameson.
\newblock A new class of high-order energy stable flux reconstruction schemes.
\newblock {\em Journal of Scientific Computing}, 47(1):50--72, 2011.

\end{thebibliography}

\end{document}